%% file: main.tex
\documentclass[journal]{new-aiaa}
\usepackage[utf8]{inputenc}
\usepackage{textcomp}

\usepackage{graphicx}
\usepackage{amsmath}
\usepackage{comment}
\usepackage[version=4]{mhchem}
\usepackage{siunitx}
\usepackage{longtable,tabularx}
\renewcommand{\arraystretch}{1.0}

\title{Direct Pseudospectral Optimal Control by Orthogonal Polynomial Integral Collocation}

\author{Thomas L. Ahrens\footnote{Graduate Research Assistant, Aerospace Engineering, 400 Bizzell St. Member AIAA. Corresponding Author, tahrens@tamu.edu.}}
\affil{Texas A\&M University, College Station, Texas, 77843}

\author{Ian M. Down\footnote{Ph.D. Candidate, Aerospace Engineering, 400 Bizzell St. Member AIAA.}}
\affil{Texas A\&M University, College Station, Texas, 77843}

\author{Manoranjan Majji\footnote{Professor, Director of Land, Air, and Space Robotics Laboratory, Aerospace Engineering, 400 Bizzell St. Associate Fellow AIAA.}}
\affil{Texas A\&M University, College Station, Texas, 77843}

\begin{document}

\maketitle

\begin{abstract}
This paper details a methodology to transcribe an optimal control problem into a nonlinear program for generation of the trajectories that optimize a given functional by approximating only the highest order derivatives of a given system's dynamics. The underlying method uses orthogonal polynomial integral collocation by which successive integrals are taken to approximate all lower order states. Hence, one set of polynomial coefficients can represent an entire coordinate's degree of freedom. Specifically, Chebyshev polynomials of the first and second kind and Legendre polynomials are used over their associated common interpolating grids derived from the bases' roots and extrema. Simple example problems compare different polynomial bases' performance to analytical solutions. The planar circular orbit raising problem is used to verify the method with solutions obtained by other pseudospectral methods in literature. Finally, a rocket landing flip maneuver problem is solved to demonstrate the ability to solve complex problems with multiple states and control variables with constraints. Simulations establish this method's performance, and reveal that the polynomial/node choice for a given problem notably affects the performance.
\end{abstract}

\input{01_intro}
\input{02_MP}
\input{03_OPIC}
\input{04_problems}
\input{05_conclusion}


\section*{Acknowledgments}
The authors thank Dr. Caleb Peck of Sandia National Laboratories for helpful discussions and for laying a framework upon which much of this research is based. Prof. John L. Junkins of Texas A\&M is also acknowledged for valuable insight and support. The authors are grateful to the partial support provided by the project entitled ``Multi-Phenomenological, Autonomous, and Understandable SDA and XDA Decision Support,'' Space University Research Initiative (SURI) Award Number FA9550-23-1-072. Principal Investigator Prof. Marcus Holzinger and Air Force Office of Scientific Research (AFOSR) Program Manager Dr. Erik P. Blasche are acknowledged. 

\bibliography{references}

\end{document}

%% file: 01_intro.tex
\section{Introduction}


\lettrine{T}{rajectory} optimization is integral to the successful use of aerospace vehicles --- particularly those whose mission criteria feature time and fuel sensitive constraints. Because of this, dynamic trajectory optimization continues to be an active area of both academic and applied research. Methods of trajectory optimization for both continuous and discrete time dynamical systems generally fall into two categories: direct methods and indirect methods. Indirect methods convert the dynamic optimization problem into a two (multi) point boundary value problem (TPBVP) using calculus of variations and Pontryagin's minimum principle. The resulting first order necessary conditions formulated in terms of the differential equations associated with the sensitivity of the cost functional with respect to the state variables guarantee solutions that are found to be local extrema in most normal problems \cite{bryson2018applied}. This TPBVP is typically solved using either a shooting or collocation scheme with the dynamic constraint variables termed costates or adjoints whose dimension mirrors that of the problem's state vector. Direct methods convert a dynamic optimization problem into one of parameter optimization which can be solved as a nonlinear program (NLP) using a variety of problem agnostic algorithms \cite{conway2010spacecraft}. Rich literature exists to describe the equivalency between the direct and indirect approaches \cite{fahroo2002direct}. Lagrange multipliers are used in the enforcement of constraints when solving a nonlinear program. It has been shown that these multipliers from the direct problem link back to the indirect problem's costates when certain direct pseudospectral methods are used \cite{gong2008connections,gong2010costate,Spada}. This then allows for direct methods to initialize indirect methods for solution refinement in the continuous domain.

Both approaches to trajectory optimization have their own unique set of advantages and challenges. Indirect methods have relatively low computational dimension and converged solutions are guaranteed to be locally optimal by the satisfaction of the first order necessary conditions. Unfortunately, solutions to problems formulated this way have small convergence basins and are hypersensitive to initial guesses of costates, making them difficult to seed. This fundamental flaw is magnified when a problem definition includes state dependent path constraints, discontinuities in the dynamics, and/or actuation constraints. Direct methods are more robust to this phenomena but can have a high computational dimension, and suffer from the curse of dimensionality in methods which discretize the time domain. Solutions obtained from solving a direct approach via conversion into an NLP are locally optimal but are ultimately approximations of continuous time solutions. In that regard, the choice of a direct or indirect approach is largely application dependent and often complementary. For a comprehensive survey on approaches to trajectory optimization, see Betts \cite{betts1998survey} and Topputo \cite{topputo2014survey}.

Among the variety of direct optimization methods, pseudospectral (PS) methods have proven to be robust, efficient, and applicable to a wide range of problems, which is why they form the foundation for several commercially available optimization software suites such as GPOPS-II \cite{patterson2014gpops}, DIDO \cite{ross2020enhancements}, SPARTAN \cite{SPARTAN}, and PSOPT \cite{PSOPT}. They exhibit spectral accuracy, or the property that truncation error from domain discretization decays as quickly as the global smoothness of the underlying solution allows \cite{CPNumAna_minimax}. PS methods in optimal control date back to the 1980s \cite{vlassenbroeck1988chebyshev}, and established themselves in a high stakes application when they were used to solve/execute a constrained zero propellant reorientation of the International Space Station in 2006 \cite{ZPM}. Since then, much research has been done in exploring new approaches to PS optimal control, connections between direct PS and indirect methods, and adaptive domain segmentation. However, the general approach of direct PS methods has largely remained the same: choose an orthogonal polynomial basis to represent a problem's states (and sometimes control), and a corresponding grid of nodes over which to discretize a problem's independent variable domain. Then leverage derivative and/or integral properties of the polynomial basis' generating functions to create a set of dynamical (collocation) constraints at each node. Lastly, use quadrature weights corresponding to the polynomial basis and grid to evaluate the problem's continuous time cost function in a discrete manner. A review of direct PS optimal control and its foundations is detailed by Ross and Karpenko \cite{ross2012review}.

The derivative approach by which discretized state and control histories are represented by a basis of polynomials with unknown coefficients is most common. The derivative properties of a particular polynomial evaluated on a corresponding grid are leveraged so that state derivatives may be computed without matrix inversion and dynamics can be collocated at each grid node. Useful grids are usually related to Legendre or Chebyshev polynomial zeros, extrema, or combinations of the two. The Legendre PS method \cite{elnagar1995pseudospectral,ross2004legendre} uses Legendre-Gauss-Lobatto (LGL) nodes, while the Gauss PS method \cite{benson2006direct,benson2005gauss} uses Legendre-Gauss (LG) nodes, and the Radau PS method \cite{garg2011direct} uses Legendre-Gauss-Radau (LGR) nodes.  All use Lagrange interpolating polynomials to arrive at an explicit differentiation matrix as a function of the Legendre polynomial basis. It is important to note that the nodes on the aforementioned Legendre polynomial grids do not have closed form, explicit representations, and must be solved with numerical methods. The LGL nodes contain both domain boundary points, while the LGR nodes contain one, and the LG nodes fail to contain any. Without nodes explicitly on the domain boundary, boundary constraints must be enforced through interpolation and/or quadrature. The choice of grid with respect to Legendre polynomials also has consequences on the resulting solution's accuracy and costate mapping back to the indirect problem. It has been shown that the Gauss and Radau PS methods yield accuracy several orders of magnitude better than the Legendre PS method, while also providing better approximations to the corresponding indirect solution's costates \cite{garg2017overview}. All three methods can be formulated in an integral form as noted in \cite{garg2017overview} where the dynamics (instead of the states) are represented by Lagrange interpolating polynomials in terms of a Legendre basis evaluated on the appropriate grid. The polynomial's integral relationships are then used to force collocated equivalence through the integration of the dynamics --- though this is only explicitly shown in \cite{benson2005gauss}. Chebyshev PS methods have been developed in the same sense, over a Chebyshev-Gauss (CG) grid \cite{vlassenbroeck1988chebyshev}, and over a Chebyshev-Gauss-Lobatto (CGL) grid using single differentiation \cite{fahroo2002direct}, and using cell-averaging and higher order differentiation \cite{el1995chebyshev}. These methods ultimately use Chebyshev polynomials of the first kind as their interpolating polynomial basis. Chebyshev polynomials of the second kind employed in direct PS optimal control are much less common, and are not readily found in literature related to aerospace applications. An example of Chebyshev polynomials of the second kind used in solving fractional optimal control problems (OCPs) is presented in \cite{nemati2016spectral}. More recently, integral properties of Gegenbauer polynomials are used to transcribe the OCP leveraging successive integration properties \cite{gegenBVP,gegenDirect}. This method showed superior performance index results for several example problems when compared to the other methods mentioned above. After an OCP is transcribed into an NLP using one of the direct PS methods above, the problem's free variables become the state values at each point on the grid along with the control (and relevant times). This is referred to as a \textit{nodal} approach. When differentiation is used to collocate the dynamics, these methods generally require the dynamics be cast into first order state space form such that successive derivatives of a single discretized state are avoided.
 
This paper introduces a new direct PS transcription of the OCP by way of orthogonal polynomial integral collocation (OPIC) using Chebyshev and Legendre polynomials. OPIC extends the principles of integral Chebyshev collocation (ICC) \cite{Peck:2023} to CGL nodes, Chebyshev polynomials of the second kind, and the LG and LGL nodes of Legendre polynomials, with several important differences that are discussed in the following sections. The underlying method is also well suited to solve initial value problems (IVP) and boundary value problems (BVP), as detailed in \cite{Peck:2023}. This method is most similar to the approach in \cite{gegenDirect}, but deviates in its general applicability to all orthogonal polynomial bases, giving a user the ability to compare different polynomials whose bases may span the solution space better than others. In contrast to the mathematical operation of differentiation, integration is widely known as a smoothing procedure that reduces computational error. Using integration as the foundation for enforcing collocation constraints allows for state dynamics to be cast into their highest order form so that only one set of polynomial coefficients is used to represent all derivative levels of a particular state. This in turn leads to a significantly smaller computational dimension in the transcribed NLP. In the proposed transcription process, the highest order derivative of each state (or degree of freedom) is represented by a series of Chebyshev or Legendre polynomials, and integral properties associated with the polynomial basis and an appropriate grid are leveraged to compute approximations of lower level state derivatives. The choice of Chebyshev or Legendre grid with respect to the polynomial basis is critical to mitigate the increase in order of approximation inherent in the integration properties. In contrast to the \textit{nodal} approaches referenced above, the proposed OPIC approach is \textit{modal}: the free variables related to the states are the unknown polynomial coefficients used to approximate the highest level state derivatives. The control component remains the same: discretized values at each node along the appropriate grid. Control remains nodal to ensure bang-bang or discontinuous control solutions can be effectively computed.

The paper is organized as follows. Section~\ref{sec:math_prelim} covers the mathematics necessary for the implementation of this direct PS method, including orthogonal function approximation, Chebyshev and Legendre polynomials, and Gaussian quadrature. Section~\ref{sec:OPIC} describes the OPIC method and direct transcription of the OCP using OPIC. In Section \ref{sec:numerical_examples}, simple example problems illustrate the accuracy and performance of this method compared to existing analytical solutions, and computation time and accuracy performance metrics amongst the different polynomial bases and nodes are also presented. Then, several difficult aerospace problems are solved to demonstrate the capabilities of the method. Lastly, key findings are reiterated to the reader in Section~\ref{sec:conclusion}. Though an important tool in robust PS implementation, domain segmentation, or in other words, the choice of local versus global collocation \cite{huntington2008comparison,darby2011hp} is not examined in this work. And while the NLP that solves the collocation problem does produce a set of Lagrange multipliers, their mapping back to the indirect problem is not covered here.

%% file: 02_MP.tex
\section{Mathematical Preliminaries} \label{sec:math_prelim}
This section serves as a mathematical and notational primer for the following sections. First, the basics of orthogonal polynomials and their use in function approximation is reviewed. Then the orthogonal polynomials of interest in this paper are presented, namely Chebyshev polynomials of the first and second kind and Legendre polynomials. Finally, the quadrature weights associated with the different interpolation nodes used in the analysis are discussed.

\subsection{Chebyshev and Legendre Polynomials} \label{sec:ChebyshevLegendrePoly}
In the context of numerical methods, particularly collocation methods, using specific families of orthogonal polynomials can significantly enhance the accuracy and efficiency of solutions. For finite-domain problems, polynomials defined with orthogonality conditions over a similarly finite interval (e.g., $[-1,1]$) are used. Additionally, certain integral properties are necessary for the formulation presented in the following section. Thus, the polynomials used in this analysis are Chebyshev and Legendre polynomials. This is not an exhaustive list of polynomials that benefit from the following formulation, but rather contains the polynomials most popular in PS literature.

\subsubsection*{Chebyshev Polynomials}
Chebyshev polynomials are a class of orthogonal polynomials highly regarded for numerical approximation due to several of their inherent properties. One of their key advantages is their near-optimal minimization of the approximation error, which stems from cosine spaced interpolation nodes, avoiding the typical pitfalls of polynomial interpolation, such as the well-known Runge phenomenon \cite{boyd2001chebyshev}. The most commonly used subset are Chebyshev polynomials of the first kind, but Chebyshev polynomials of the second kind have also been shown to be useful, particularly for numerical integration \cite{MasonChebyPoly}. Both polynomials are defined over the domain $\tau \in [-1, 1]$. Given that the $n\tth$ order polynomial has exactly $n$ zeros (or roots) and $n+1$ extrema over the Chebyshev domain $\tau \in [-1,1]$, these are obvious choices for the nodes points $\boldsymbol{\tau}$.

Chebyshev polynomials of the first kind (CP1K) can be defined in many ways, but given their close relation to the sine and cosine trigonometric functions, they are written as such in Eq.~\eqref{eq:CP1K}.  These polynomials can also be defined with a recursive generating function. The first two polynomials are defined in Eqs.~\eqref{eq:CP1K_0} and \eqref{eq:CP1K_1}, and the recurrence relation for CP1K is given in Eq.~\eqref{eq:CP1K_n}. CP1K satisfy the continuous orthogonality condition with the weighting function $w(\tau) = (1 - \tau^2)^{-\frac{1}{2}}$.
\begin{align}
    T_n(\tau) &= \cos (n \arccos \tau), \qquad \tau \in [-1,1] \label{eq:CP1K} \\
    T_0(\tau) &= 1 \label{eq:CP1K_0} \\
    T_1(\tau) &= \tau \label{eq:CP1K_1} \\
    T_{n+1}(\tau) &= 2\tau T_{n}(\tau) - T_{n-1}(\tau), \qquad n \geq 1 \label{eq:CP1K_n}
\end{align}

As mentioned above, two sets of nodes are of interest for CP1K. The zeros, or Chebyshev-Gauss (CG) nodes, are defined in Eq.~\eqref{eq:CP1K_zeros} as the roots of the $(n+1)\tth$ order Chebyshev polynomial. The extrema, or Chebyshev-Guass-Lobatto (CGL) nodes, are defined in Eq.~\eqref{eq:CP1K_ext} as the interior extrema and endpoints $\pm 1$ of the $n\tth$ order Chebyshev. The $(n+1)\tth$ order polynomial is used for the CG nodes such that the final column of the matrix ${}^{n}\Phi(\boldsymbol{\tau})$ is not zero, since the $n\tth$ Chebyshev polynomial evaluated at its own roots yields zeros, and so that the number of CG nodes is equivalent to the number of CGL nodes. The CG nodes have been the subject of recent work \cite{Peck:2023, SciTech}, and have been shown to be optimal interpolation nodes, as supported by Clenshaw's minimax theorem \cite{CPNumAna_minimax}. The CGL nodes have also been used in various methods for solving IVPs and BVPs \cite{Bai}.
\begin{align}
    \tau_k^\text{CG} &= \cos \left( \frac{(k + \frac{1}{2}) \pi}{n+1} \right), \qquad k = [0, 1, \ldots, n] \label{eq:CP1K_zeros} \\
    \tau_k^\text{CGL} &= \cos \left( \frac{k \pi}{n} \right), \qquad k = [0, 1, \ldots, n] \label{eq:CP1K_ext}
\end{align}
When the matrix ${}^{n}\Phi(\boldsymbol{\tau})$ uses the CG nodes as the grid points $\boldsymbol{\tau}$, the columns are discretely orthogonal with respect to the unit weighting function $w(\tau_k) = 1$. With the same weighting function, the CGL nodes produce discretely orthogonal columns when the first and last terms, corresponding to the endpoints $\pm 1$ for $k=[0,n]$, are halved, denoted by the double prime notation $\Sigma {}^{\prime\prime}$ in Eq.~\eqref{eq:CP1K_disc_orthog} \cite{MasonChebyPoly}.
\begin{align}
    \sum_{k=1}^{n+1} T_i(\tau_k^\text{CG}) T_j(\tau_k^\text{CG}) &= \sum_{k=0}^{n} {}^{\prime\prime} \ T_i(\tau_k^\text{CGL}) T_j(\tau_k^\text{CGL}) = \begin{cases} 0, & i \neq j; \: \{i, j\} \leq n \\ \frac{n}{2}, & i = j \neq 0; \: \{i, j\} \leq n \\ n, & i = j = 0 \end{cases} \label{eq:CP1K_disc_orthog}
\end{align}

Integrals relating the polynomials to their adjacent polynomials are leveraged in the collocation formulation of the following sections, so these integrals are presented here. The first two indefinite integrals of CP1K are easily observed and given in Eqs.~\eqref{eq:cp1k_int0} and \eqref{eq:cp1k_int1}, while the subsequent integrals follow the form of Eq.~\eqref{eq:cp1k_intn} \cite{boyd2001chebyshev, RussianIntegrals}.
\begin{align}
    \int T_0 (\tau) \dif{\tau} &= T_1(\tau) \label{eq:cp1k_int0} \\
    \int T_1 (\tau) \dif{\tau} &= \frac{1}{4} \left( T_2(\tau) + T_0(\tau) \right) \label{eq:cp1k_int1} \\
    \int T_n (\tau) \dif{\tau} &= \frac{1}{2} \left( \frac{T_{n+1}(\tau)}{n+1} - \frac{T_{n-1}(\tau)}{n-1} \right), \qquad n \geq 2 \label{eq:cp1k_intn}
\end{align}

Chebyshev polynomials of the second kind (CP2K) are defined trigonometrically in Eq.~\eqref{eq:CP2K}, and with the first two polynomials in Eqs.~\eqref{eq:CP2K_0} and \eqref{eq:CP2K_1}, follow the recurrence generating function of Eq.~\eqref{eq:CP2K_n}. These polynomials satisfy continuous orthogonality with the weighting function $w(\tau) = (1 - \tau^2)^{\frac{1}{2}}$.
\begin{align}
    U_n(\tau) &= \frac{\sin ((n+1) \arccos \tau)}{\sin (\arccos \tau)}, \qquad \tau \in [-1,1] \label{eq:CP2K} \\
    U_0 (\tau) &= 1 \label{eq:CP2K_0} \\
    U_1 (\tau) &= 2 \tau \label{eq:CP2K_1} \\
    U_{n+1}(\tau) &= 2\tau U_{n}(\tau) - U_{n}(\tau), \qquad n \geq 1 \label{eq:CP2K_n}
\end{align}
CP2K can be related to CP1K in several ways, including that of Eq.~\eqref{eq:CP1K_CP2K_relation}. Additionally, CP1K can be differentiated to obtain CP2K using Eq.~\eqref{eq:CP1K_deriv_CP2K}.
\begin{align}
    U_n(\tau)  &= U_{n-2}(\tau) + 2 T_n(\tau), \qquad n \geq 2 \label{eq:CP1K_CP2K_relation} \\
    T_n^{\prime}(\tau) &= n U_{n-1}(\tau) \label{eq:CP1K_deriv_CP2K}
\end{align}

The zeros of the the $(n+1)\tth$ order polynomial defined in Eq.~\eqref{eq:CP2K_zeros} are the nodes of interest. While the extrema of CP1K are clearly discernible, the extrema of CP2K are not. However, the extrema of CP2K with respect to the their continuous weighting function $w(\tau) = (1 - \tau^2)^{\frac{1}{2}}$ can be found explicitly \cite{MasonChebyPoly}, but these nodes are not considered here. Note that the zeros of the $(n+1)\tth$ order CP2K are equivalent to the interior extrema of the $n \tth$ order CP1K.
\begin{align}
    \tau_k &= \cos \left( \frac{k \pi}{n+2} \right), \qquad k = [1, 2, \ldots, n+1] \label{eq:CP2K_zeros}
\end{align}
Using the CP2K nodes of Eq.~\eqref{eq:CP2K_zeros} in the matrix ${}^{n}\Phi(\boldsymbol{\tau})$, the columns are discretely orthogonal with respect to the weighting function $w(\tau_k) = (1 - \tau_k^2)$ when the first and last terms are halved, as shown in Eq.~\eqref{eq:CP2K_disc_orthog} \cite{MasonChebyPoly}.
\begin{align}
    \sum_{k=0}^{n} {}^{\prime\prime} (1-\tau_k^2) U_i(\tau_k) U_j(\tau_k) &= \begin{cases} 0, & i \neq j; \: \{i, j\} \leq n-1 \\ \frac{n}{2}, & i = j \neq 0; \: 0 \leq \{i, j\} \leq n-1 \end{cases} \label{eq:CP2K_disc_orthog}
\end{align}

As with CP1K, the first two integrals of CP2K are easily observed in Eqs.~\eqref{eq:cp2k_int0} and \eqref{eq:cp2k_int1}, and the subsequent integrals are given in Eq.~\eqref{eq:cp2k_intn} \cite{RussianIntegrals}.
\begin{align}
    \int U_0 (\tau) \dif{\tau} &= \frac{1}{2} U_1(\tau) \label{eq:cp2k_int0} \\
    \int U_1 (\tau) \dif{\tau} &= \frac{1}{4} \left( U_2(\tau) + U_0(\tau) \right) \label{eq:cp2k_int1} \\
    \int U_n (\tau) \dif{\tau} &= \frac{1}{2(n+1)} \left( U_{n+1} - U_{n-1} \right), \qquad n \geq 2 \label{eq:cp2k_intn}
\end{align}

\subsubsection*{Legendre Polynomials}
Another set of orthogonal polynomials commonly used in numerical analysis are Legendre polynomials $P_n(\tau)$. These polynomials satisfy continuous orthogonality condition with the weighting function $w(\tau)=1$. The first two Legendre polynomials are defined in Eqs.~\eqref{eq:Leg0} and \eqref{eq:Leg1}, and the recurrence relation is given in Eq.~\eqref{eq:Legn}.
\begin{align}
    P_0(\tau) &= 1 \label{eq:Leg0} \\
    P_1(\tau) &= \tau \label{eq:Leg1} \\
    P_{n+1}(\tau) &= \frac{2n+1}{n+1} \tau P_n(\tau) - \frac{n}{n+1} P_{n-1}(\tau), \qquad n \geq 1 \label{eq:Legn}
\end{align}

The nodes of interest for the Legendre polynomials are Legendre-Gauss (LG) and Legendre-Gauss-Lobatto (LGL) nodes. The LG nodes are defined as the zeros of the $(n+1)\tth$~Legendre polynomial, and the LGL nodes are the zeros of the derivative of the $n\tth$~Legendre polynomial, $P_n^\prime$, plus the endpoints $\pm 1$. Unlike the nodes used for the Chebyshev polynomials, however, these nodes cannot be defined explicitly, and must be found numerically. Although upper and lower bounds on the position of the roots have been refined over the last century \cite{SzegoBounds}, many fast algorithms have been developed in recent years to find these nodes \cite{FORTRAN_quad, HaleTownsend, Bogaert}. The integrals of the first two Legendre polynomials are observed in Eqs.~\eqref{eq:legint0} and \eqref{eq:legint1}, and the subsequent integrals are given in Eq.~\eqref{eq:legintn} \cite{RussianIntegrals}.
\begin{align}
    \int P_0(\tau) \dif{\tau} &= P_1(\tau) \label{eq:legint0} \\
    \int P_1(\tau) \dif{\tau} &= \frac{1}{3} P_2(\tau) - \frac{1}{6} P_0(\tau) \label{eq:legint1} \\
    \int P_n(\tau) \dif{\tau} &= \frac{1}{2n+1} \left( P_{n+1}(\tau) - P_{n-1}(\tau) \right), \qquad n \geq 2 \label{eq:legintn}
\end{align}

\subsection{Quadrature Weights} \label{sec:quadweights}
An often necessary component of numerical analysis is numerical integration or quadrature, where a definite integral is approximated by some method. In the context of optimal control, this can relate to approximating a continuous integral cost function. Perhaps the best known method is Gaussian quadrature, wherein an integral is approximated as a linear combination of weights and the integrand evaluated at specific nodes. Gaussian quadrature was first expanded to orthogonal polynomials by Jacobi, where the definite integral is converted to the domain $[-1,1]$ and the rule takes the form of Eq.~\ref{eq:GaussLegQuad}, where $n$ is the number of samples, $w_k$ are the quadrature weights, and $\tau_k$ are specific node points, typically the roots of a particular polynomial \cite{GaussQuad}.
\begin{align} \label{eq:GaussLegQuad}
    \int_{-1}^{1} f(\tau) \dif{\tau} &\approx \sum_{k=1}^{n} w_k f(\tau_k)
\end{align}
In this form, this well-studied method is known as Gauss-Legendre quadrature, and there are many modern algorithms which efficiently compute these nodes and weights \cite{GolubWelsch}. When using the LG and LGL nodes, both of which have $n+1$ points, the quadrature weights are given in Eqs.~\eqref{eq:LGWeight} and \eqref{eq:LGLWeight}, respectively \cite{AbramowitzStegun, SpectralMethodsFluidDynamics, MethNumInt}.
\begin{align}
    w_k^\text{LG} &= \frac{2}{(1-\tau_k^2) \left[ P_{n+1}^\prime (\tau_k) \right]^2}, \qquad k = [0, 1, \ldots, n] \label{eq:LGWeight} \\
    w_k^\text{LGL} &= \frac{2}{n(n+1)} \frac{1}{\left[ P_n (\tau_k) \right]^2}, \qquad k = [0, 1, \ldots, n] \label{eq:LGLWeight}
\end{align}

This quadrature rule was extended to CP1K via a change of variables and discrete cosine transform in a method known as Clenshaw-Curtis or Fej\'er quadrature \cite{clenshaw1960method, Fejer}, and was subsequently used to derive the CP2K weights \cite{MethNumInt}. The quadrature weights using the $n+1$ CG, CGL, and CP2K nodes are given in Eqs.~\eqref{eq:cp1k_cg_weights}, \eqref{eq:cp1k_cgl_weights}, and \eqref{eq:cp2k_weights}, respectively. For all weight formulae, $\theta_k$ are related to the interpolation nodes as $\theta_k = \cos^{-1} \tau_k$, and $\lfloor \cdot \rfloor$ denotes the floor function or the integer component of a positive real number.
\begin{align}
    w_k^\text{CG} &= \frac{2}{n+1} \left\{ 1 - \sum_{i=0}^{\lfloor (n+1)/2 \rfloor} \frac{\cos ( 2i \theta_k )}{4i^2 - 1} \right\}, \qquad k = [0, 1, \ldots, n] \label{eq:cp1k_cg_weights}
\end{align}
\vspace{-3.5em}
\begin{subequations} \label{eq:cp1k_cgl_weights}
\begin{align}
    w_k^\text{CGL} &= \frac{2}{n} \left\{ 1 - \sum_{i=0}^{\lfloor n/2 \rfloor} \frac{\cos(2i \theta_k)}{4i^2 - 1} \right\}, \qquad k = [1, \ldots, n-1] \qquad \:\:\: \\
    w_0^\text{CGL} &= w_n^\text{CGL} = \begin{cases} \frac{1}{n^2}, & \text{if $n$ is even} \\ \frac{1}{(n-1)(n+1)}, & \text{if $n$ is odd} \end{cases}
\end{align}
\end{subequations}
\vspace{-3.5em}
\begin{align}
    w_k^\text{CP2K} &= \frac{4 \sin \theta_k}{n+2} \sum_{i=0}^{\lfloor (n+2)/2 \rfloor} \frac{\sin((2i - 1) \theta_k)}{2i - 1}, \qquad k = [0, 1, \ldots, n]  \label{eq:cp2k_weights}
\end{align}

In order to use these quadrature rules, the integral must first be converted from the time domain $t \in [t_0, t_f]$ to the computational domain $\tau \in [-1,1]$ such that it is in the form of Eq.~\eqref{eq:GaussLegQuad}. This is accomplished with the linear, isomorphic mapping of Eq.~\eqref{eq:time_transform}, where $\Delta t = t_f-t_0$. Time derivatives of any level taken in either domain can then be related by Eq. \ref{eq:dtdtau}. Thus, integrals are converted between domains as Eq.~\eqref{eq:Int_trans}
\begin{align}
    \tau = \frac{2}{\Delta t} t - \frac{t_f+t_0}{\Delta t} \qquad &\Leftrightarrow \qquad t = \frac{\Delta t}{2}\tau + \frac{t_f + t_0}{2} \label{eq:time_transform} \\
    \left( \frac{\dif{\tau}}{\dif{t}} \right)^n = \left( \frac{2}{\Delta t} \right)^n \qquad &\Leftrightarrow \qquad \left( \frac{\dif{t}}{\dif{\tau}} \right)^n = \left( \frac{\Delta t}{2} \right)^n \label{eq:dtdtau} \\
    \int_{t_0}^{t_f} f(t) \dif{t} &= \frac{\Delta t}{2} \int_{-1}^{1} f(\tau) \dif{\tau} \label{eq:Int_trans}
\end{align}

%% file: 03_OPIC.tex
\section{Orthogonal Polynomial Integral Collocation} \label{sec:OPIC}
Orthogonal Polynomial Integral Collocation (OPIC) is a novel technique in which the dynamics (i.e., highest order derivative) are approximated by a series of orthogonal polynomials. Integral relationships relating the adjacent polynomials are then leveraged to integrate the same polynomial series until reaching the desired derivative level, which is typically that of the state-space. Unlike other integral methods as in \cite{Elgindy}, no additional coefficients are needed for each integration: a single set describes every derivative level for a given coordinate. The integral properties of OPIC are exact, and the matrix formulation leads to a much simpler solution when many integrations are required. A similar method with several minor, but important, formulation differences was derived in \cite{Peck:2023, SciTech} with CP1K using the CG nodes. In this section, a new generalized formulation is developed, which works not only for CG nodes, but also the CGL nodes, Chebyshev polynomials of the second kind, and Legendre polynomials using the LG and LGL nodes. A generic optimal control problem is then presented, and the direct OPIC transcription is given.

\subsection{Dynamics Approximation and Integral Formulation}
To form a solution using OPIC, the highest order derivative $y^{(q)}$ at the point $\tau \in [-1,1]$ is represented by an $n\tth$ order series of CP1K, CP2K, or Legendre polynomials, written as
\begin{align}
    y^{(q)} (\tau) &= \sum_{i=0}^{n} \phi_i (\tau) \alpha_i = {}^n\boldsymbol{\phi}(\tau) \bm{\alpha} \label{eq:y^q} \\
    {}^n\boldsymbol{\phi}(\tau) &\in \mathbb{R}^{1 \times (n+1)}, \quad \bm{\alpha} \in \mathbb{R}^{(n+1) \times 1} \nonumber
\end{align}

Now that the dynamics are represented, the next highest derivative is given by integrating the dynamics series representation from $-1$ to $\tau$.
\begin{align}
    \int_{-1}^{\tau} y^{(q)} (\tau) \dif{\tau} &= \int_{-1}^{\tau} {}^n\boldsymbol{\phi}(\tau) \bm{\alpha} \dif{\tau} \\
    y^{(q-1)}(\tau) - y^{(q-1)}(-1) &= \left[ \int_{-1}^{\tau} {}^n\boldsymbol{\phi}(\tau) \dif{\tau} \right] \bm{\alpha}
\end{align}
\begin{table}[th!]
    \caption{\label{tab:int_weights} Integration weights for Chebyshev and Legendre polynomials}
    \centering
    \renewcommand{\arraystretch}{1.5} 
    \begin{tabular}{lccccc}
        \hline \hline
        Polynomial, $\phi(\tau)$ & $b_0^+$ & $b_1^+$ & $b_1^-$ & $b_i^+$ & $b_i^-$ \\ \hline
        CP1K, $T_n(\tau)$ & $1$ & $\frac{1}{4}$ & $\frac{1}{4}$ & $\frac{1}{2(i+1)}$ & $-\frac{1}{2(i-1)}$ \\
        CP2K, $U_n(\tau)$ & $\frac{1}{2}$ & $\frac{1}{4}$ & $\frac{1}{4}$ & $\frac{1}{2(i+1)}$ & $-\frac{1}{2(i+1)}$ \\
        Legendre, $P_n(\tau)$ & $1$ & $\frac{1}{3}$ & $-\frac{1}{6}$ & $\frac{1}{2i+1}$ & $-\frac{1}{2i+1}$ \\
        \hline \hline
    \end{tabular}
\end{table}
The constant coefficients $\bm{\alpha}$ can be pulled outside so the integral can be taken term by term for each polynomial in the series. Following the orthogonal polynomial recursive relations, the integrals of these polynomials can be expressed in terms of the adjacent polynomials, where the zeroeth, first and $i\tth$ integrals are written as
\begin{subequations}
\begin{align}
    \int \phi_0 (\tau) \dif{\tau} &= b_0^+ \phi_1(\tau) + C_0 \\
    \int \phi_1 (\tau) \dif{\tau} &= b_1^+ \phi_{2} (\tau) + b_1^- \phi_{0} (\tau) + C_1 \\
    \int \phi_i (\tau) \dif{\tau} &= b_i^+ \phi_{i+1} (\tau) + b_i^- \phi_{i-1} (\tau) + C_i, \qquad i \geq 2
\end{align}
\end{subequations}
$b_0^+$, $b_1^+$, $b_1^-$, $b_i^+$, and $b_i^-$ are the integration weights for the adjacent polynomials, with the subscript denoting the order of the polynomials being integrated and the superscript denoting whether the weight corresponds to the higher or lower order polynomial resulting from integration.  The specific weights for each polynomial given in Sec.~\ref{sec:ChebyshevLegendrePoly} are compiled in Table~\ref{tab:int_weights}. These integral formulae are then applied to the entire polynomial series as a definite integral from $-1$ to $\tau$.
\begin{align} 
    \int_{-1}^{\tau} {}^{n}\boldsymbol{\phi}(\tau) \dif{\tau} &= \int_{-1}^{\tau} \begin{bmatrix} {\phi}_0(\tau) & {\phi}_1(\tau) & \cdots & {\phi}_n(\tau) \end{bmatrix} \dif{\tau} \nonumber \\
    &= \left[ b_0^+ \phi_1(\tau) + \left( b_1^+ \phi_{2}(\tau) + b_1^- \phi_{0}(\tau) \right) + \ldots + \left( b_n^+ \phi_{n+1}(\tau) + b_n^- \phi_{n-1}(\tau) \right) \right]_{-1}^{\tau} \nonumber \\
    &= \left[ b_1^- \phi_0(\tau) + \left( b_0^+ + b_2^- \right) \phi_{1}(\tau) + \left( b_1^+ + b_3^- \right) \phi_{2}(\tau) + \ldots \right. \nonumber \\
    &\quad \left. + \left( b_{n-2}^+ + b_{n}^- \right) \phi_{n-1}(\tau) + b_{n-1}^+ \phi_{n}(\tau) + b_{n}^+ \phi_{n+1}(\tau) \right]_{-1}^{\tau} \nonumber \\
    \int_{-1}^{\tau} {}^{n}\boldsymbol{\phi}(\tau) \dif{\tau} &= \left[ {}^{n+1} \boldsymbol{\phi}(\tau) - {}^{n+1} \boldsymbol{\phi} (-1) \right] \mathcal{B}_{1} \label{eq:int_scheme}
\end{align}
The row vectors ${}^{n+1}\boldsymbol{\phi}(\tau_k)$ and ${}^{n+1}\boldsymbol{\phi}(-1)$ are defined by appending an $(n+1)$ order polynomial to the end of the row vector, and $\mathcal{B}_{1}$ is an integration weight matrix composed of entries only on the super- and sub-diagonals given as
\begin{align}
    \mathcal{B}_{1} &= \begin{bmatrix} 0 & b_1^- & 0 & \cdots & 0 & 0 & 0 \\ b_0^+ & 0 & b_2^- & \ddots & 0 & 0 & 0 \\ 0 & b_1^+ & 0 & \ddots & 0 & 0 & 0 \\ 0 & 0 & b_2^+ & \ddots & 0 & 0 & 0 \\ \vdots & \vdots & \vdots & \ddots & \vdots & \vdots & \vdots \\ 0 & 0 & 0 & \ddots & b_{n-2}^+ & 0 & b_n^- \\ 0 & 0 & 0 & \ddots & 0 & b_{n-1}^+ & 0 \\ 0 & 0 & 0 & \cdots & 0 & 0 & b_n^+ \end{bmatrix} \in \mathbb{R}^{(n+2) \times (n+1)} \label{eq:B1}
\end{align}
Thus, the next lowest derivative is given as
\begin{align}
    y^{(q-1)}(\tau) &= \left[ {}^{n+1}\boldsymbol{\phi}(\tau) \mathcal{B}_{1} - \bm{v}_1 \right] \bm{\alpha} + y^{(q-1)}(-1) \label{eq:yq1} \\
    \bm{v}_1 &= {}^{n+1}\boldsymbol{\phi}(-1) \mathcal{B}_{1} \in \mathbb{R}^{1 \times (n+1)} \label{eq:v1}
\end{align}
where $\bm{v}_1$ is the evaluated lower integration bound which is unchanged for all $\tau$.

To yield the second lowest derivative, Eq.~\eqref{eq:yq1} is integrated as
\begin{align}
    \int_{-1}^{\tau} y^{(q-1)}(\tau) \dif{\tau} &= \int_{-1}^{\tau} \left( \left[ {}^{n+1}\boldsymbol{\phi}(\tau) \mathcal{B}_{1} - \bm{v}_1 \right] \bm{\alpha} + y^{(q-1)}(-1) \right) \dif{\tau} \\
    y^{(q-2)}(\tau) - y^{(q-2)}(-1) &= \left[ \int_{-1}^{\tau} \left( {}^{n+1}\boldsymbol{\phi}(\tau) \mathcal{B}_{1} - \bm{v}_1 \right) \dif{\tau} \right] \bm{\alpha} + (\tau + 1) y^{(q-1)}(-1)
\end{align}
The remaining integral is given as
\begin{align}
    \int_{-1}^{\tau} \left( {}^{n+1}\boldsymbol{\phi}(\tau) \mathcal{B}_{1} - \bm{v}_1 \right) \dif{\tau} &= \left( \int_{-1}^{\tau} {}^{n+1}\boldsymbol{\phi}(\tau) \dif{\tau} \right) \mathcal{B}_{1} - (\tau+1) \bm{v}_1  \\
    &= \left( \left[ {}^{n+2}\boldsymbol{\phi}(\tau) - {}^{n+2}\boldsymbol{\phi}(-1) \right] \mathcal{B}_{2} \right) \mathcal{B}_{1} - (\tau+1) \bm{v}_1 
\end{align}
where $\mathcal{B}_{2}$ is the integration weight matrix arising from integrating ${}^{n+1}\boldsymbol{\phi}(\tau)$. Since the integration weights for the $0\tth$ through $n\tth$ polynomials remained unchanged, only one additional row and column are needed to define this matrix succinctly in a block matrix form.
\begin{align}
    \mathcal{B}_{2} &= \begin{bmatrix} \mathcal{B}_{1} & \rvline & \begin{matrix} \bm{0}_{n \times 1} \\ b_{n+1}^- \\ 0 \end{matrix} \\ \hline \bm{0}_{1 \times (n+1)} & \rvline & b_{n+1}^+ \end{bmatrix} \label{eq:B2}
\end{align}
Thus, the second highest derivative is
\begin{align}
    y^{(q-2)}(\tau) &= \left[ {}^{n+2}\boldsymbol{\phi}(\tau) \mathcal{B}_{2} \mathcal{B}_{1} - \bm{v}_{2} - (\tau + 1) \bm{v}_{1} \right] \bm{\alpha} + (\tau + 1) y^{(q-1)}(-1) + y^{(q-2)}(-1) \label{eq:yq2} \\
    \bm{v}_{2} &= {}^{n+2}\boldsymbol{\phi}(-1) \mathcal{B}_{2} \mathcal{B}_{1} \in \mathbb{R}^{1 \times (n+1)} \label{eq:v2}
\end{align}

This integration process can be continued to generate an equation for the $(q-m)\tth$ derivative, which is given in the following form.
\begin{align}
    y^{(q-m)}(\tau) &= \left[ {}^{n+m}\boldsymbol{\phi}(\tau) \left( \prod_{j=m}^{1} B_j \right) - \sum_{j=1}^{m} p_{m-j}(\tau) \bm{v}_{j} \right] \bm{\alpha} + \sum_{j=1}^{m} p_{m-j}(\tau) y^{(q-j)}(-1) \label{eq:yqm} \\
    \bm{v}_{j} &= {}^{n+j}\boldsymbol{\phi}(-1) \left( \prod_{h=j}^{1} B_h \right) \in \mathbb{R}^{1 \times (n+1)} \label{eq:vj} \\
    p_{m-j}(\tau) &= \frac{1}{(m-j)!} \sum_{l=0}^{m-j} \binom{m-j}{l} \tau^{m-j-l} \in \mathbb{R}^{1 \times 1} \label{eq:pmj} \\
    B_j &= \begin{bmatrix} \mathcal{B}_{j-1} & \rvline & \begin{matrix} \bm{0}_{(n+j-2) \times 1} \\ b_{n+j-1}^- \\ 0 \end{matrix} \\ \hline \bm{0}_{1 \times (n+j-1)} & \rvline & b_{n+j-1}^+ \end{bmatrix} \label{eq:Bj}
\end{align}
The $\bm{v}_{j}$ vectors are constant lower integration bounds, $B_j$ are the integration weight matrices, $y^{(q-j)}(-1)$ are the initial conditions of the different derivative levels, and $p_{m-j}(\tau)$ are the terms arising from integrating a constant from $-1$ to $\tau$. For convenience, the first few $p_{m-j}(\tau)$ terms are given below.
\begin{subequations}
\begin{align}
    p_0(\tau) &= 1 \\
    p_1(\tau) &= (\tau + 1) \\
    p_2(\tau) &= \frac{1}{2} (\tau^2 + 2\tau + 1) \\
    p_3(\tau) &= \frac{1}{6} (\tau^3 + 3\tau^2 + 3\tau + 1)
\end{align}
\end{subequations}

\subsection{Collocation Procedure}
The integral relations derived above allow for any level below the dynamics to represented using a single set of coefficients. Thus, dynamics of any order can be represented with an orthogonal polynomial series and subsequently integrated until reaching the desired derivative level. While the formula of Eq.~\eqref{eq:yqm} presents a global solution that holds for any $\tau \in [-1,1]$, the coefficients must first be found to employ it. To achieve this, a collocation approach is applied. This procedure is simple: replace all instances of $\tau$ with $\bm{\tau}$, the set of nodes corresponding to the chosen orthogonal basis functions. The formulation is then written in vector-matrix form as
\begin{align}
    \bm{y}^{(q-m)}(\bm{\tau}) &= \left[ {}^{n+m}\Phi(\bm{\tau}) \left( \prod_{j=m}^{1} B_j \right) - \sum_{j=1}^{m} \bm{p}_{m-j}(\bm{\tau}) \bm{v}_{j} \right] \bm{\alpha} + \sum_{j=1}^{m} \bm{p}_{m-j}(\bm{\tau}) y^{(q-j)}(-1) \label{eq:yqm_colloc}
\end{align}
where
\begin{align}
    {}^{n+m}\Phi(\bm{\tau}) &= \begin{bmatrix} \phi_0(\tau_0) & \phi_1(\tau_0) & \cdots & \phi_n(\tau_0) & \phi_{n+1}(\tau_0) & \cdots & \phi_{n+m}(\tau_0) \\ \phi_0(\tau_1) & \phi_1(\tau_1) & \cdots & \phi_n(\tau_1) & \phi_{n+1}(\tau_1) & \cdots & \phi_{n+m}(\tau_1) \\ \vdots & \vdots & \cdots & \vdots & \vdots & \cdots & \vdots \\ \phi_0(\tau_n) & \phi_1(\tau_n) & \cdots & \phi_n(\tau_n) & \phi_{n+1}(\tau_n) & \cdots & \phi_{n+m}(\tau_n) \end{bmatrix} \in \mathbb{R}^{(n+1) \times (n+m+1)} \\[2em]
    \bm{p}_{m-j}(\bm{\tau}) \bm{v}_{j} &= \begin{bmatrix} p_{m-j}(\tau_0) \bm{v}_{j} \\ p_{m-j}(\tau_1) \bm{v}_{j} \\ \vdots \\ p_{m-j}(\tau_n) \bm{v}_{j} \end{bmatrix} \in \mathbb{R}^{(n+1) \times (n+1)} \\
    \bm{p}_{m-j}(\bm{\tau}) y^{(q-j)}(-1) &= \begin{bmatrix} p_{m-j}(\tau_0) y^{(q-j)}(-1) \\ p_{m-j}(\tau_1) y^{(q-j)}(-1) \\ \vdots \\ p_{m-j}(\tau_n) y^{(q-j)}(-1) \end{bmatrix} \in \mathbb{R}^{(n+1) \times 1}
\end{align}

This procedure is efficient since all of the matrices and vectors defined above are unchanged once the approximation order is chosen. For high order systems, OPIC does not require that the system be rewritten in state-space form. Instead of having $m$ sets of coefficients, only one set is needed to fully describe all derivative levels. Additionally, integration acts as a smoothing operation that can be helpful numerically. OPIC can be used in solving IVPs and BVPs, wherein BVPs are converted into a TPBVP. To solve IVPs, the coefficients $\bm{\alpha}$ can be found through a single matrix inversion for a linear system. Nonlinear systems require iterative methods such as successive substitution or Newton's method. In the following section, the direct transcription for solving optimal control problems is presented.

\subsection{Direct Transcription of the Optimal Control Problem}
The general, continuous time optimal control problem expressed in Bolza form is now converted into an NLP with direct OPIC (DOPIC). This process is termed "direct transcription." The Bolza form is typically shown with dynamics in state space form. To make clear the advantage of DOPIC as a PS optimal control method, the transcription that follows deals with a second order, scalar dynamical system --- without loss of generality. Consider the following Bolza problem with cost function $J$, second order dynamics $f$, and boundary constraints $\bm{\psi}$
\begin{align}
    \min \qquad J &= \phi(y_0, \dot{y}_0, t_0, y_f, \dot{y}_f, t_f) + \int_{t_0}^{t_f} g(y, \dot{y}, u, t) \dif{t} \label{eq:bolza} \\
\begin{split}
    \text{s.t.} \qquad \ddot{y} &= f(y, \dot{y}, u, t) \\
    \bm{\psi}_0 (y_0, \dot{y}_0, t_0) &= \bm{0} \\
    \bm{\psi}_f (y_f, \dot{y}_f, t_f) &= \bm{0}
\end{split}
\end{align}
where $u(t)$ is the control time history minimizing $J$ over the time horizon.

The transcription begins by transforming the time domain $t \in [t_0, t_f]$ of Eq.~\ref{eq:bolza} to the computational domain $\tau \in [-1,1]$ over which the above orthogonal polynomials operate. This is accomplished using the mapping of Eq.~\eqref{eq:time_transform}, and time derivatives are taken as Eq. \ref{eq:dtdtau}. This transcription process is shown with global collocation such that one polynomial basis is used to approximate the entire domain. Domain segmentation can be introduced to solve problems with interior point constraints or discrete events, but that is beyond the scope of this method's introduction.

Next, the highest order dynamics (in this case, second order dynamics) are approximated in Eq.~\eqref{eq:ypp_approx} with an $n\tth$ order orthogonal polynomial basis at the corresponding grid nodes $\bm{\tau}$, where $\bm{\alpha}$ are the set of $n+1$ unknown modal coefficients, and $(\cdot)^{\prime}$ denotes $\frac{\dif{}}{\dif{\tau}}$. Using the previously derived integral formulae for the first and second integral, the lower level derivative and state are then given in Eqs.~\eqref{eq:ypstep} and \eqref{eq:ystep}, respectively. Here, initial conditions are assumed to be known with $y(-1) = y(t_0)$ and $y^{\prime}(-1) = \frac{\Delta t}{2}\dot{y}(t_0)$.
\begin{align}
    \bm{y}^{\prime\prime}(\bm{\tau}) &= {}^{n}\Phi(\bm{\tau}) \bm{\alpha} \label{eq:ypp_approx} \\
    \bm{y}^{\prime}(\bm{\tau}) &= \left[ {}^{n+1}\Phi(\bm{\tau}) \mathcal{B}_{1} - \bm{p}_0(\bm{\tau}) \bm{v}_1 \right] \bm{\alpha} + \bm{p}_0(\bm{\tau}) y^{\prime}(-1) \label{eq:ypstep} \\
    \bm{y}(\bm{\tau}) &= \left[ {}^{n+2}\Phi(\bm{\tau}) \mathcal{B}_{2} \mathcal{B}_{1} - \bm{p}_0(\bm{\tau}) \bm{v}_{2} - \bm{p}_1(\bm{\tau}) \bm{v}_{1} \right] \bm{\alpha} + \bm{p}_1(\bm{\tau}) y^{\prime}(-1) + \bm{p}_0(\bm{\tau}) y(-1) \label{eq:ystep}
\end{align}

Unlike the state representation, the controls are represented nodally such that each node $\tau_k$ has a control variable $u_k$ placed there exactly. There is no assumption on the control behavior behavior between nodes since the values at the nodes are all that is needed for the collocation procedure. With expressions for all derivative levels, the dynamic collocation constraint $\xi_k$ of Eq.~\eqref{eq:dyn_con_k} is formed at the $k\tth$ node $\tau_k$ which contains the unknown modal coefficients and nodal controls. Collecting the dynamic collocation constraints at the full set of nodes $\bm{\tau}$, this can be written succinctly as the constraint vector $\Xi$ of Eq.~\eqref{eq:all_dyn_con}.
\begin{align}
    \xi_k &= {}^{n}\boldsymbol{\phi} (\tau_k) \bm{\alpha} - \left( \frac{\Delta t}{2} \right)^2 f \left( y_k, y_k^{\prime}, u_k, \tau_k \right) = 0 \label{eq:dyn_con_k} \\
    \Xi &= \begin{bmatrix} \xi_0 & \xi_1 & \cdots & \xi_n \end{bmatrix}^\top \in \mathbb{R}^{n+1 \times 1} = {}^{n}\Phi(\bm{\tau}) \bm{\alpha} - \left( \frac{\Delta t}{2} \right)^2 \bm{f} \left( \bm{y}, \bm{y}^{\prime}, \bm{u}, \bm{\tau} \right) = \bm{0} \label{eq:all_dyn_con}
\end{align}
Note that although this equality occurs in the problem's transformed computational domain, the dynamics themselves are not changed. That is, $\dot{\bm{y}}$ must be transformed into $\bm{y}^{\prime}$ before direct substitution into the dynamics. Lastly, boundary constraints $\bm{\psi}_0$ and $\bm{\psi}_f$ are converted to the computational domain and are included in the discretized formulation. Grids which do not contain nodes on the boundaries, namely CG, CP2K, and LG, use the interpolated solution to evaluate the endpoints $\tau = \pm 1$, which is readily accessible from the coefficients. Though not explicitly shown in this section, nonlinear state/control dependent path equality and inequality constraints are included in the same manner.

To minimize the continuous cost functional in Eq. \ref{eq:bolza}, the integral is discretized into a summation over the corresponding computational grid for evaluation after the limits are transformed. The approximation then takes the form 
\begin{align}
    \frac{\Delta t}{2} \int_{-1}^{1} g (y, y^{\prime}, u, \tau) \dif{\tau} &\approx \frac{\Delta t}{2} \sum_{k=0}^{n} g_k (y_k, y_k^{\prime}, u_k, \tau_k) \: w_k \\
    \Tilde{J} = \phi + \frac{\Delta t}{2} &\sum_{k=0}^{n} g_k (y_k, y_k^{\prime}, u_k, \tau_k) \: w_k
\end{align}
where $w_k$ are the quadrature weights associated with the corresponding nodes defined in Sec.~\ref{sec:quadweights}.

With the continuous time optimal control problem transcribed using DOPIC, the resulting NLP is stated simply as
\begin{align}
    \min_{\chi} \Tilde{J}(\chi) \quad \text{s.t.} \quad \Ceq = \bm{0}, \quad \Cin \leq \bm{0}
\end{align}
The free variable vector $\chi$ and nonlinear equality constraint vector $\Ceq$ are
\begin{align}
    \chi &= \begin{bmatrix} \bm{\alpha} \\ \bm{u} \\t_0 \\ t_f \end{bmatrix} \in \mathbb{R}^{(2(n+1)+4) \times 1}, \qquad \Ceq = \begin{bmatrix} \bm{\psi}(-1) \\ \Xi \\ \bm{\psi}(1) \end{bmatrix} = \bm{0} \in \mathbb{R}^{(n+5) \times 1}
\end{align}
assuming $\bm{\psi}(-1)$ and $\bm{\psi}(1)$ contain boundary constraints on the state and derivative at initial and final time, respectively. A generalized inequality constraint vector $\Cin$ captures any inequality constraints which may be placed on the state, derivatives, control, time, or other variables.

%% file: 04_problems.tex
\section{Numerical Examples} \label{sec:numerical_examples}
In this section, several classic trajectory optimizations problems are first solved using DOPIC to demonstrate the method's performance against problems with analytical solutions. A second order minimum fuel problem with control inequality constraints is solved, followed by Breakwell's problem with state inequality constraints. While these problems can be reduced to linear or quadratic programs, they are represented as NLPs to exhibit the robustness of the method. Then, two aerospace specific problems are solved with DOPIC to show the method's ability to solve complex problems. First, a constant thrust, planar circular orbit raising problem is solved for the minimum time and maximum radius cases. Finally, a simplified version of the Starship landing flip maneuver problem is used to demonstrate the ability to a multi-coordinate problem with constraints on both the states and controls. All computations are are run on an Apple M1 Max\textsuperscript{\textregistered} CPU operating at 3.228 GHz using MATLAB\textsuperscript{\textregistered} 2024a's {\fontfamily{qcr}\selectfont fmincon} function with finite-differenced Jacobians. All problems are solved using the SQP method except the rocket landing problem, which uses the interior-point method.

\subsection{Second Order Minimum Fuel with Control Inequality}
The first problem is a one degree-of-freedom, second-order minimum fuel problem with control inequalities. The problem has a closed form analytical solution derived from an indirect method, and includes a ``bang-zero-bang'' control structure with two switching times \cite{BrysonHo}. The problem aims to find $u(t)$ to minimize
\begin{align}
    J &= \int_{0}^{t_f} |u(t)| \dif{t} \label{eq:P1cost} \\
\begin{split} \label{eq:P1dyn}
    \text{s.t.} \qquad \ddot{x}(t) &= u(t), \qquad -1 \leq u \leq 1 \\
    x(0) &= x_0, \qquad \dot{x}(0) = \dot{x}_0 \\
    x(t_f) &= 0, \qquad \dot{x}(t_f) = 0
\end{split}
\end{align}
As in \cite{bryson2018applied, BrysonHo}, the case considered is one with the initial conditions and final time of Eq.~\eqref{eq:SO_ICs}, where $t_{f,{\min}}$ is the minimum time to take the system from the initial state $(x_0, \dot{x}_0)$ to the final state $(0,0)$ with the control inequality constraint $-1 \leq u \leq 1$.
\begin{align} \label{eq:SO_ICs}
    \dot{x}_0 &\geq 0, \qquad x_0 \geq -\frac{1}{2} \dot{x}_0^2, \qquad t_f \geq t_{f,{\min}} = \dot{x}_0 + \sqrt{4x_0+2\dot{x}_0^2}
\end{align}

To solve this problem using DOPIC, one $n\tth$ order series is used to approximate the dynamics of Eq.~\eqref{eq:P1dyn} with controls at the interpolation nodes. Thus, the free variable vector $\chi$, equality constraint vector $\Ceq$, and inequality constraint vector $\Cin$ are defined as
\begin{align} \label{eq:P1vectors}
    \chi &= \begin{bmatrix} \bm{\alpha} \\ \bm{u} \end{bmatrix} \in \mathbb{R}^{2(n+1) \times 1}, \quad \Ceq = \begin{bmatrix} x(-1) - x_0 \\ x^\prime (-1) - \left( \frac{t_f}{2} \right) \dot{x}_0 \\ \Xi \\ x(1) - x_f \\ x^\prime (1) - \left( \frac{t_f}{2} \right) \dot{x}_f \end{bmatrix} = \bm{0} \in \mathbb{R}^{(n+5) \times 1}, \quad \Cin = \begin{bmatrix} \bm{u} - 1 \\ -1 - \bm{u} \end{bmatrix} \leq \bm{0} \in \mathbb{R}^{2(n+1) \times 1}
\end{align}

\begin{figure}[hbt!]
    \centering
    \includegraphics[width=.95\textwidth]{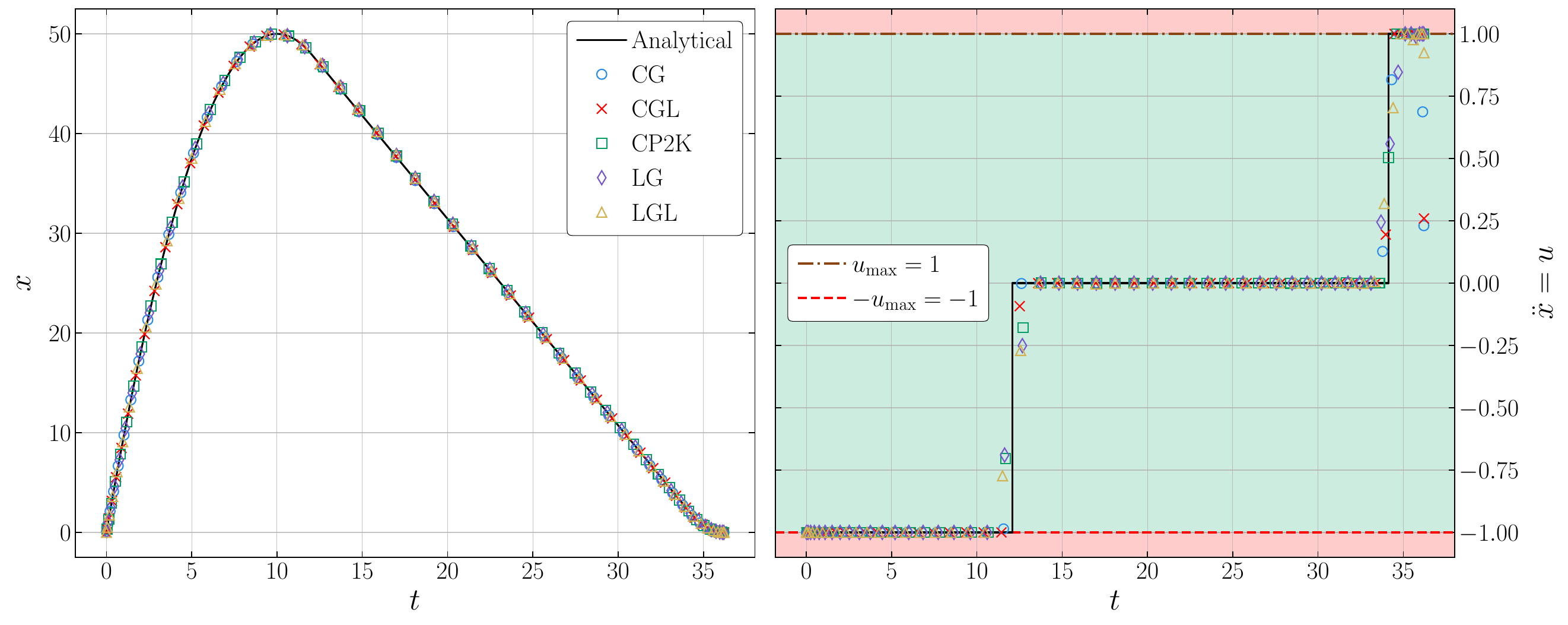}
    \caption{State and control history for the second order minimum control problem solved with $n=50$. Constrained regions are shown in red, while the unconstrained region is shown in green.}
    \label{fig:SOMinControl_states}
\end{figure}

\begin{figure}[hbt!]
    \centering
    \includegraphics[width=.95\textwidth]{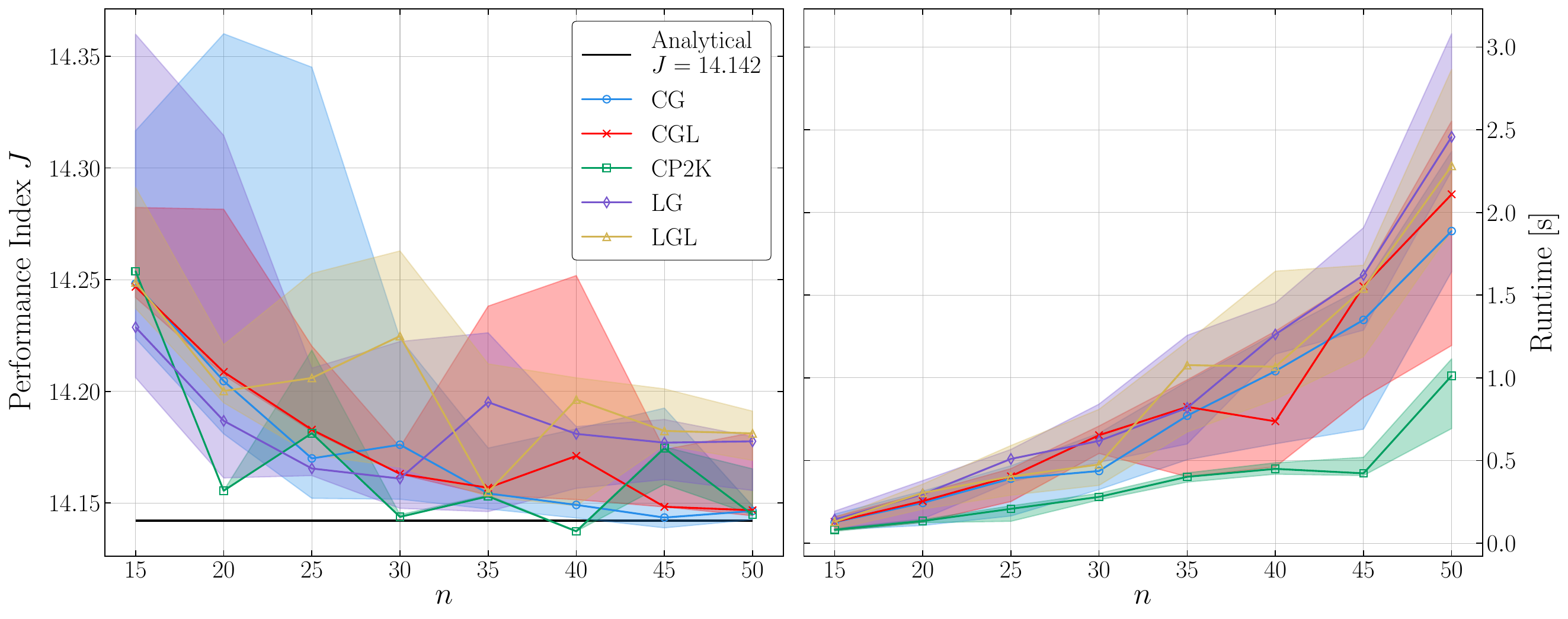}
    \caption{Performance index and runtime quartiles for the second order minimum control problem solved with approximation order from $n=15$ to $n=50$.}
    \label{fig:SOMinControl_perf}
\end{figure}

The state and control history for all polynomial and node sets with an order $n=50$ solution are shown in Fig.~\ref{fig:SOMinControl_states}. Both histories follow the analytical solution closely, although the approximation struggles at the discontinuities in control at the switching times. However, this is expected behavior for this problem's discontinuous solution. With increasing approximation order $n$, the DOPIC solution better captures the discontinuities.

The performance of the DOPIC solution was analyzed by 100 runs for eight values of $n$, increasing in increments of 5 from 15 to 50, and the results were compared to the analytical solution. The system values used are $x_0 = 0$, $\dot{x}_0 = 10$, and $t_f = 1.5 \left( t_f \right)_{\min}$. The free variable vector $\chi$ for each run was initialized by a guess of uniformly random elements between zero and one. The performance index and computation time results are shown in Fig.~\ref{fig:SOMinControl_perf}. All approximation orders converged for the provided random initial conditions, demonstrating the robustness of the DOPIC solution method. As expected, computational runtime increases and the performance index approaches the true cost as the approximation order increases. The lower runtime for CP2K suggests that these basis functions may be better suited for this problem specifically.

\subsection{Breakwell Problem}
The Breakwell problem is a minimum-energy, second order problem with state variable inequality constraints \cite{BrysonHo}. The problem is to find $u(t)$ for $t \in [0,1]$ to minimize
\begin{align}
    J &= \frac{1}{2} \int_{0}^{1} u^2 (t) \dif{t} \label{eq:P2cost} \\
\begin{split} \label{eq:P2dyn}
    \text{s.t.} \qquad \ddot{x}(t) &= u(t) \\
    x(0) &= x(1) = 0 \\
    \dot{x}(0) &= \dot{x}(1) = 1 \\
    x(t) &\leq l
\end{split}
\end{align}
The Breakwell problem has an analytical solution which can be broken up into three parts: an unconstrained problem ($l \geq \sfrac{1}{4}$) with parabolic solutions, a constrained problem ($\sfrac{1}{6} \leq l \leq \sfrac{1}{4}$) with osculating cubic arc solutions, and a constrained problem ($0 < l \leq \sfrac{1}{6}$) with bifurcated cubic arc solutions \cite{BrysonHo}. Thus, any DOPIC solution must be able to approximate solutions for all three regimes.

As with the previous problem, one $n\tth$ order orthogonal polynomial series is used to approximate the dynamics of Eq.~\eqref{eq:P2dyn} with controls at the interpolation nodes. The free variable vector $\chi$, equality constraint vector $\Ceq$, and inequality constraint vector $\Cin$ are given as
\begin{align} \label{eq:P2vectors}
    \chi &= \begin{bmatrix} \bm{\alpha} \\ \bm{u} \end{bmatrix} \in \mathbb{R}^{2(n+1) \times 1}, \quad \Ceq = \begin{bmatrix} x(-1) - x_0 \\ x^\prime (-1) - \left( \frac{t_f}{2} \right) \dot{x}_0 \\ \Xi \\ x(1) - x_f \\ x^\prime (1) - \left( \frac{t_f}{2} \right) \dot{x}_f \end{bmatrix} = \bm{0} \in \mathbb{R}^{(n+5) \times 1}, \quad \Cin = \begin{bmatrix} \bm{x} (\bm{\tau}) - l \end{bmatrix} \leq \bm{0} \in \mathbb{R}^{(n+1) \times 1}
\end{align}

\begin{figure}[b!]
    \centering
    \includegraphics[width=.95\textwidth]{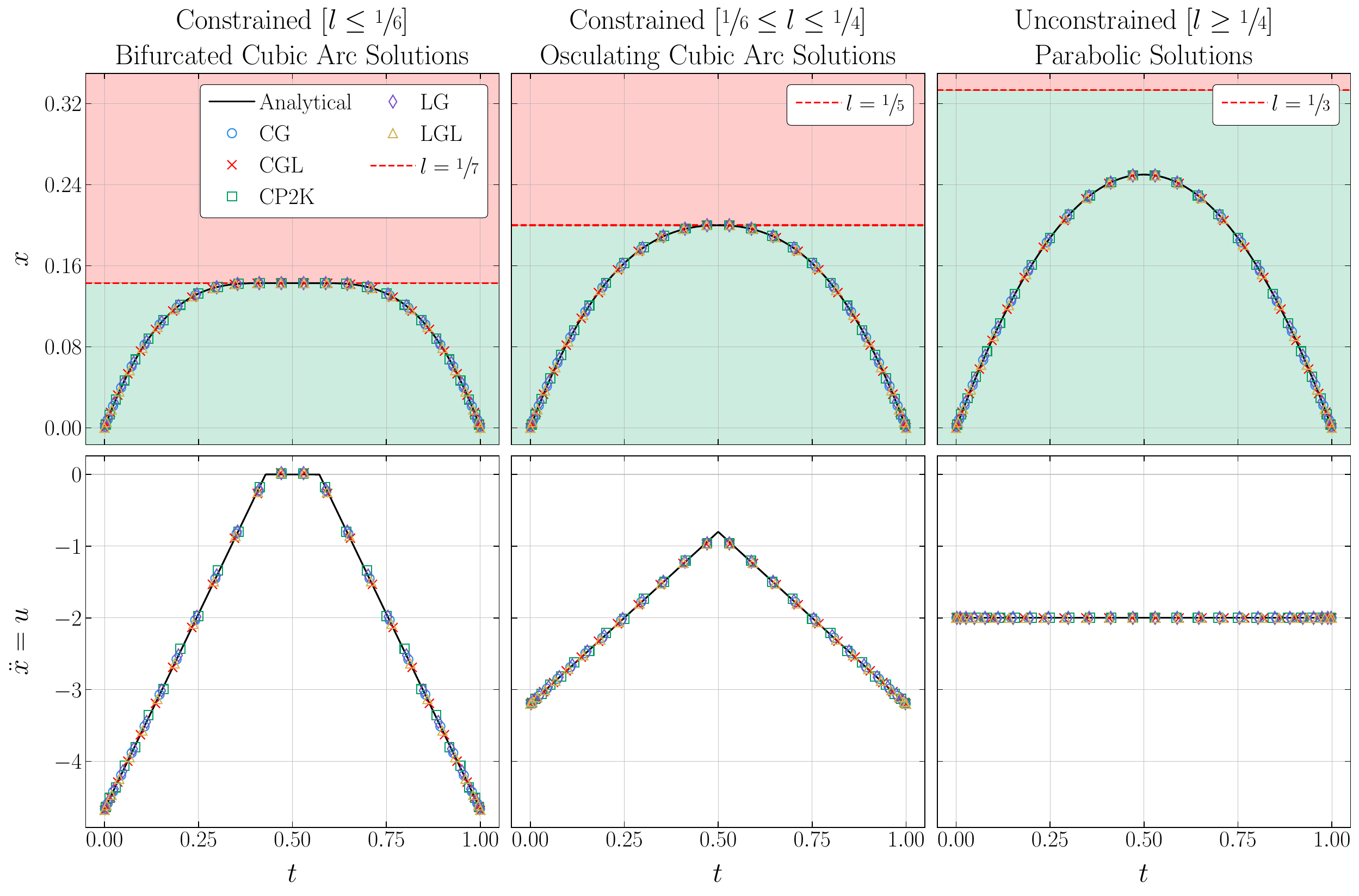}
    \caption{Breakwell problem solved for $l=[\sfrac{1}{7}, \sfrac{1}{5}, \sfrac{1}{3}]$ with $n=25$. Constrained and unconstrained regions are shown in red and green, respectively.}
    \label{fig:Breakwell_states}
\end{figure}

To test all analytical solution regimes, DOPIC solutions were also calculated for different values of $l$. The state and control histories for all polynomial and node sets with an order $n=25$ solution are shown for $l=[\sfrac{1}{7}, \sfrac{1}{5}, \sfrac{1}{3}]$ in Fig.~\ref{fig:Breakwell_states}. Both the state and control histories follow the analytical solution closely, and the DOPIC solution captures the different solution regimes well.

\begin{figure}[b!]
    \centering
    \includegraphics[width=.95\textwidth]{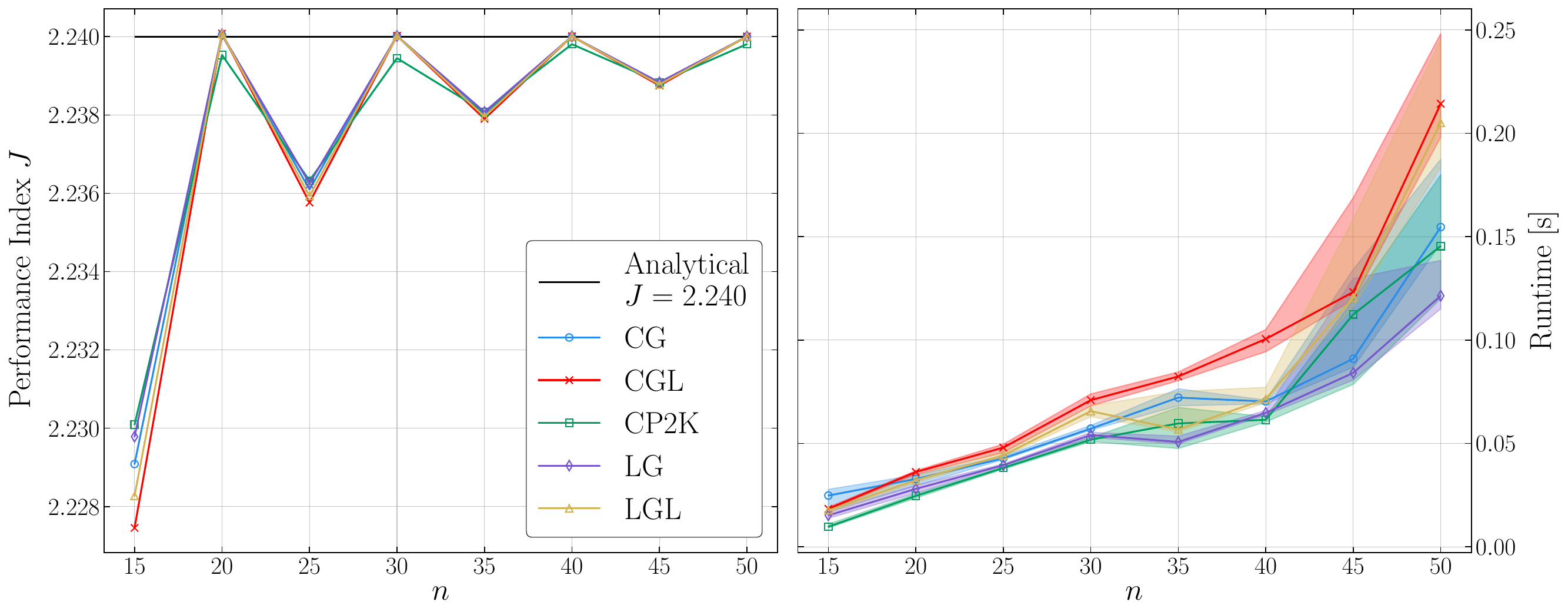}
    \caption{Performance index and runtime quartiles for the Breakwell problem solved with approximation order from $n=15$ to $n=50$ with $l = \sfrac{1}{5}$.}
    \label{fig:Breakwell_perf}
\end{figure}

Again, 100 trials were run for $n=15$ to $n=50$ to analyze the performance of the DOPIC solution, and the results were compared to the analytical solution. The value used for this analysis was $l = \sfrac{1}{5}$, and the free variable vector $\chi$ for each run was initialized by a guess of uniformly random elements between zero and one. All random initial condition sets converged, and the performance index and computation time results are shown in Fig.~\ref{fig:Breakwell_perf}. Unlike the previous problem, the cost is generally underestimated and increases with approximation order, rather than decreasing. This is due to the inequality constraints being placed on the state variable, not the control.

An interesting disparity is seen in this problem between the the solutions with even and odd approximation orders. Due to the symmetry of the problem, having an even approximation order will place a node directly on the constraint, while the odd approximation orders do not have a node directly on the constraint. Thus, for the odd order solutions, the constraint is inactive at all nodes, resulting in a less accurate performance index.

\subsection{Planar Circular Orbit Raising under Constant Thrust}
The problem of a constant, low-thrust orbit raising spacecraft between two circular planar orbits is a classical OCP \cite{elnagar1995pseudospectral, elnagar1998pseudospectral, BettsJGCD}. Two cases are considered: a minimum time Earth to Mars circular orbit transfer, and a transfer from a circular Earth orbit to the largest possible circular orbit over a fixed time interval. The cost functions and boundary conditions are different for the two cases, but both follow the same dynamics,
\begin{align}
\begin{split} \label{eq:P3dyn}
    \ddot{r}(t) &= \frac{v_t^2(t)}{r(t)} - \frac{\mu}{r^2(t)} + \frac{u_0 \sin \varphi(t)}{m_0 + \dot{m}t} \\
    \dot{v}_t(t) &= - \frac{\dot{r}(t) v_t(t)}{r(t)} + \frac{u_0 \cos \varphi(t)}{m_0 + \dot{m}t}
\end{split}
\end{align}
where $r(t)$ is the radial distance between the spacecraft and the Sun, $\dot{r}(t)$ is the radial velocity, $v_t(t)$ is the tangential velocity, and $\varphi(t)$ is the control angle as measured from the local horizon. The problem is non-dimensionalized for the gravitational parameter $\mu$, with constant low-thrust magnitude $u_0$, initial mass $m_0$, and constant propellant mass flow rate $\dot{m}$ as
\begin{align} \label{eq:non_dim_vals}
    \mu = 1, \qquad u_0 = 0.1405, \qquad m_0 = 1, \qquad \dot{m} = -0.07487
\end{align}
where the normalized time-unit (TU) is 58.18 days. The polar angle needed to fully describe the spacecraft's position in the plane is independent of this problem when the tangential velocity $v_t$ is used, and can be recovered after a solution has been found by noting the relationship
\begin{align} \label{eq:polar_recovery}
    \dot{\theta}(t) &= \frac{v_t(t)}{r(t)}
\end{align}
With the prescribed initial angle $\theta(0) = 0$, Eq.~\eqref{eq:polar_recovery} can be integrated using OPIC to get the time history of the polar angle $\theta(t)$.

To solve this problem using DOPIC, two $n\tth$ order polynomial series are used to approximate the dynamics with controls corresponding to the thrusting angle at the interpolation nodes. The first series approximates the second order dynamics of the radial coordinate, while the second series approximates the first order dynamics of the tangential coordinate.


\subsubsection{Minimum Time Orbit Transfer}
The minimum time orbit transfer problem between circular Earth and Mars orbits has the performance index
\begin{align} \label{eq:P3_mintime_cost}
    J &= \int_{0}^{t_f} 1 \dif{t} = t_f
\end{align}
At the start and end of the maneuver, the spacecraft has initial and final radial position and velocity corresponding to the circular orbits of Earth and Mars. Thus, the boundary conditions are
\begin{equation}  \label{eq:P3_mintime_BCs}
\begin{aligned}
    r(0) &= 1, \qquad &&\dot{r}(0) \: = 0, \qquad &&v_t(0) \: = 1 \\
    r(t_f) &= 1.525, \qquad &&\dot{r}(t_f) = 0, \qquad &&v_t(t_f) = 0.8098
\end{aligned}
\end{equation}
There is no exact known analytical solution for this problem, although it has been studied extensively with various PS methods \cite{elnagar1995pseudospectral, ross2012review, betts1998survey}.

A constraint on the thrusting angle $\varphi \in [0, 2\pi]$ is used, as well as a physical constraint to keep the final time non-negative, $t_f > 0$. Finally, a thrust attitude control rate constraint $\bm{c}_{\dot{\varphi}_{\max}} \in \mathbb{R}^{n \times 1}$ is implemented, which is calculated between two adjacent nodes $\tau_k$ and $\tau_j$ as
\begin{align} \label{eq:att_rate_con}
    c_{\dot{\varphi}_{\max,jk}} = \frac{|\varphi_k - \varphi_j|}{\frac{t_f}{2} \left( \tau_k - \tau_j \right)} - \dot{\varphi}_{\max} \leq 0 \quad \forall \{j, k\} \in [1, \ldots, n] \quad \text{s.t.} \quad k-j=1
\end{align}
where $\dot{\varphi}_{\max} = 400^\circ/\text{ TU}$. There are $n$ of these constraints since the attitude rate is evaluated between nodes, not at them. Thus, the free variable vector $\chi$, equality constraint vector $\Ceq$, and inequality constraint vector $\Cin$ are defined as
\begin{align} \label{eq:P3_mintime_vectors}
    \chi &= \begin{bmatrix} \bm{\alpha}_r \\ \bm{\alpha}_{v_t} \\ \bm{\varphi} \\ t_f \end{bmatrix} \in \mathbb{R}^{(3(n+1)+1) \times 1}, \quad \Ceq = \begin{bmatrix} r(-1) - r_0 \\ r^\prime (-1) - \left( \frac{t_f}{2} \right) \dot{r}_0 \\ v_t(-1) - v_{t,0} \\ \Xi_r \\ \Xi_{v_t} \\ r(1) - r_f \\ r^\prime (1) - \left( \frac{t_f}{2} \right) \dot{r}_f \\ v_t(1) - v_{t,f} \end{bmatrix} = \bm{0} \in \mathbb{R}^{(2(n+1)+6) \times 1}, \quad \Cin = \begin{bmatrix} \bm{\varphi} - 2\pi \\ -\bm{\varphi} \\ \bm{c}_{\dot{\varphi}_{\max}} \\ -t_f \end{bmatrix} \leq \bm{0} \in \mathbb{R}^{3(n+1) \times 1}
\end{align}

Given the nonlinearity of the dynamics, the problem is broken up to provide better initial free variable vector guesses. The first sub-problem does not use include the thrust attitude control rate constraint of Eq.~\eqref{eq:att_rate_con}. The random initialized polynomial coefficients $\bm{\alpha}_r$ and $\bm{\alpha}_{v_t}$ are between one and two, the values of the control angle $\varphi$ increase linearly from 0 to $\pi$, and the final normalized time guess is $t_f = 3$. Once this sub-problem has converged, the attitude rate constraints are added back in and the previous solution's final free variable vector is used as an initial guess.

\begin{figure}[t!]
    \centering
    \includegraphics[width=.95\textwidth]{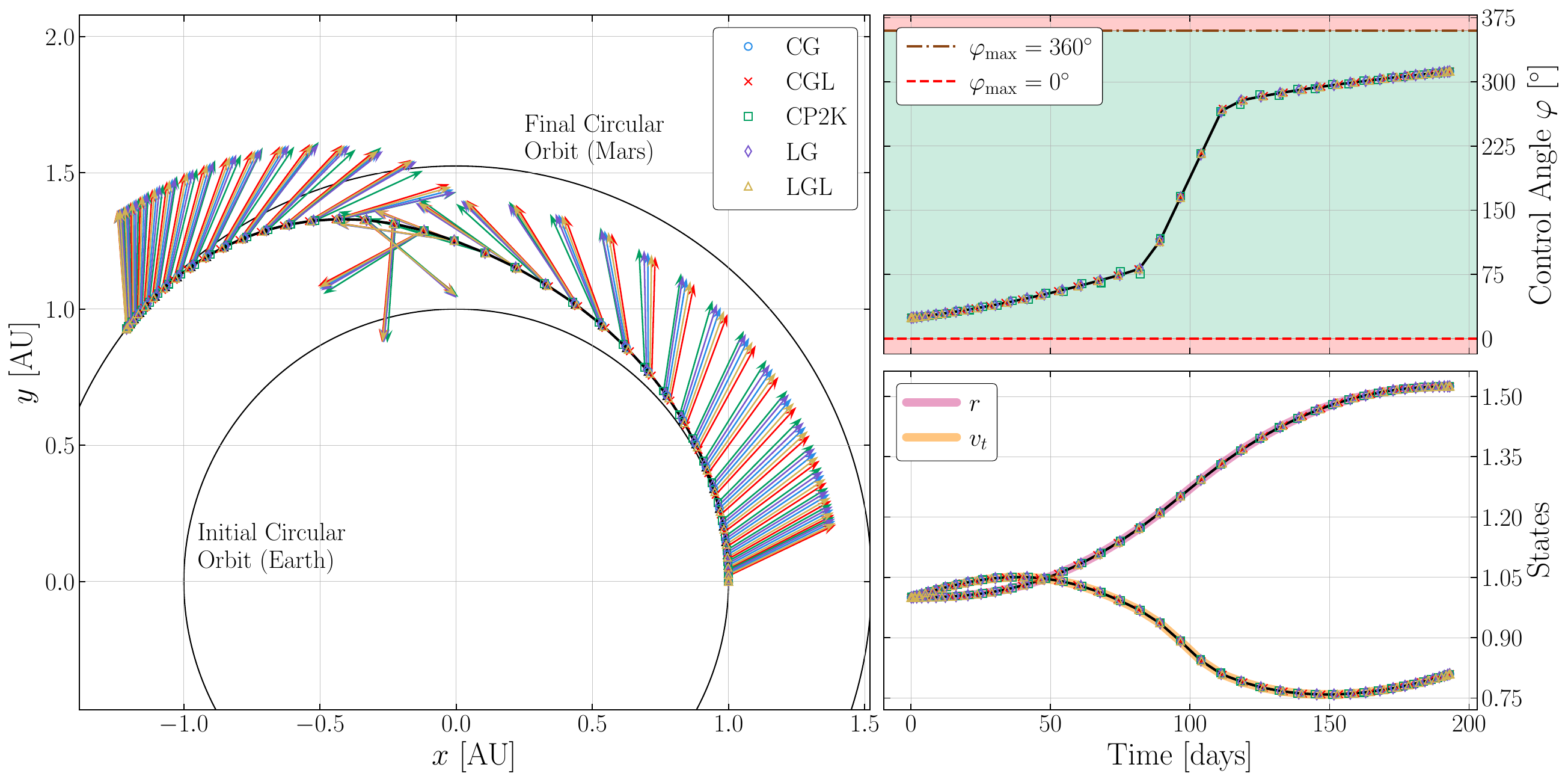}
    \caption{Minimum time planar circular orbit raising under constant thrust solved with $n=40$. Constrained regions are shown in red, while unconstrained region is shown in green.}
    \label{fig:OrbitRaisingMinTime}
\end{figure}

To analyze the performance of DOPIC on this problem, 100 trials were run for a single approximation order, $n=40$. Unlike the previous two examples, all trials did not converge due to the complexity of the problem. The number of successful runs are shown in Table~\ref{tab:EarthMars} alongside the average number of iterations, transfer time, and computational runtime for the converged trials. The orbital plane phase space and corresponding states, angles, and controls for a converged solution are shown in Fig.~\ref{fig:OrbitRaisingMinTime}. All successful runs result in similar transfer times of $t_f \simeq 3.321$, but CP2K exhibit a worse convergence percentage. The remaining polynomials and nodes have similar success rates, but the CGL nodes show higher runtimes than the others.

\begin{table}[htb!]
    \caption{\label{tab:EarthMars} DOPIC Minimum Time Orbit Transfer Performance Results Averaged Over 100 Runs for $n=40$}
    \centering
    \begin{tabular}{lcccc}
        \hline \hline
        \mc{\textbf{Nodes}} & \mc{\textbf{Successful Runs}} & \mc{\textbf{Iterations}} & \mc{\textbf{Transfer Time $t_f$ [TU]}} & \mc{\textbf{Runtime [s]}} \\ \hline
        CG   & 89  & 202 & 3.3210547 & 1.114 \\
        CGL  & 88  & 246 & 3.3210441 & 1.452 \\
        CP2K & 73  & 170 & 3.3203661 & 1.000 \\
        LG   & 90  & 180 & 3.3210516 & 0.998 \\
        LGL  & 92  & 228 & 3.3210453 & 1.147 \\
        \hline \hline
    \end{tabular}
\end{table}

\subsubsection{Maximum Radius Orbit Transfer}
The orbit transfer problem between a circular Earth orbit and a maximum radius circular orbit has the performance index
\begin{align} \label{eq:P4_maxrad_cost}
    J &= -r(t_f)
\end{align}
At the start and end of the maneuver, the spacecraft has initial radial position and velocity corresponding to the circular orbit of Earth, and at the end of the maneuver, the spacecraft is required to have a velocity corresponding to a circular orbit. Thus, the boundary conditions are
\begin{equation} \label{eq:P4_maxrad_BCs}
\begin{aligned}
    r(0) &= 1, \qquad &&\dot{r}(0) = 0, \qquad &&v_t(0) = 1 \\
    \dot{r}(t_f) &= 0, \qquad &&v_t(t_f) = \sqrt{1 / r(t_f)}
\end{aligned}
\end{equation}

The dynamics approximation remained unchanged from the minimum time problem, as do the thrusting angle and attitude control rate constraints. Since the minimum time and maximum radius problems should result in the same solution for a continuous thrust problem if the time in the two problems is the same, the final time is set to $t_f = 3.32$. The free variable vector $\chi$, equality constraint vector $\Ceq$, and inequality constraint vector $\Cin$ are now given as
\begin{align} \label{eq:P4_mintime_vectors}
    \chi &= \begin{bmatrix} \bm{\alpha}_r \\ \bm{\alpha}_{v_t} \\ \bm{\varphi} \end{bmatrix} \in \mathbb{R}^{3(n+1) \times 1}, \quad \Ceq = \begin{bmatrix} r(-1) - r_0 \\ r^\prime (-1) - \left( \frac{t_f}{2} \right) \dot{r}_0 \\ v_t(-1) - v_{t,0} \\ \Xi_r \\ \Xi_{v_t} \\ r^\prime (1) - \left( \frac{t_f}{2} \right) \dot{r}_f \\ v_t(1) - \sqrt{1/r(t_f)} \end{bmatrix} = \bm{0} \in \mathbb{R}^{(2(n+1)+5) \times 1}, \quad \Cin = \begin{bmatrix} \bm{\varphi} - 2\pi \\ -\bm{\varphi} \\ \bm{c}_{\dot{\varphi}_{\max}} \end{bmatrix} \leq \bm{0} \in \mathbb{R}^{(3(n+1)-1) \times 1}
\end{align}

\begin{figure}[ht!]
    \centering
    \includegraphics[width=.95\textwidth]{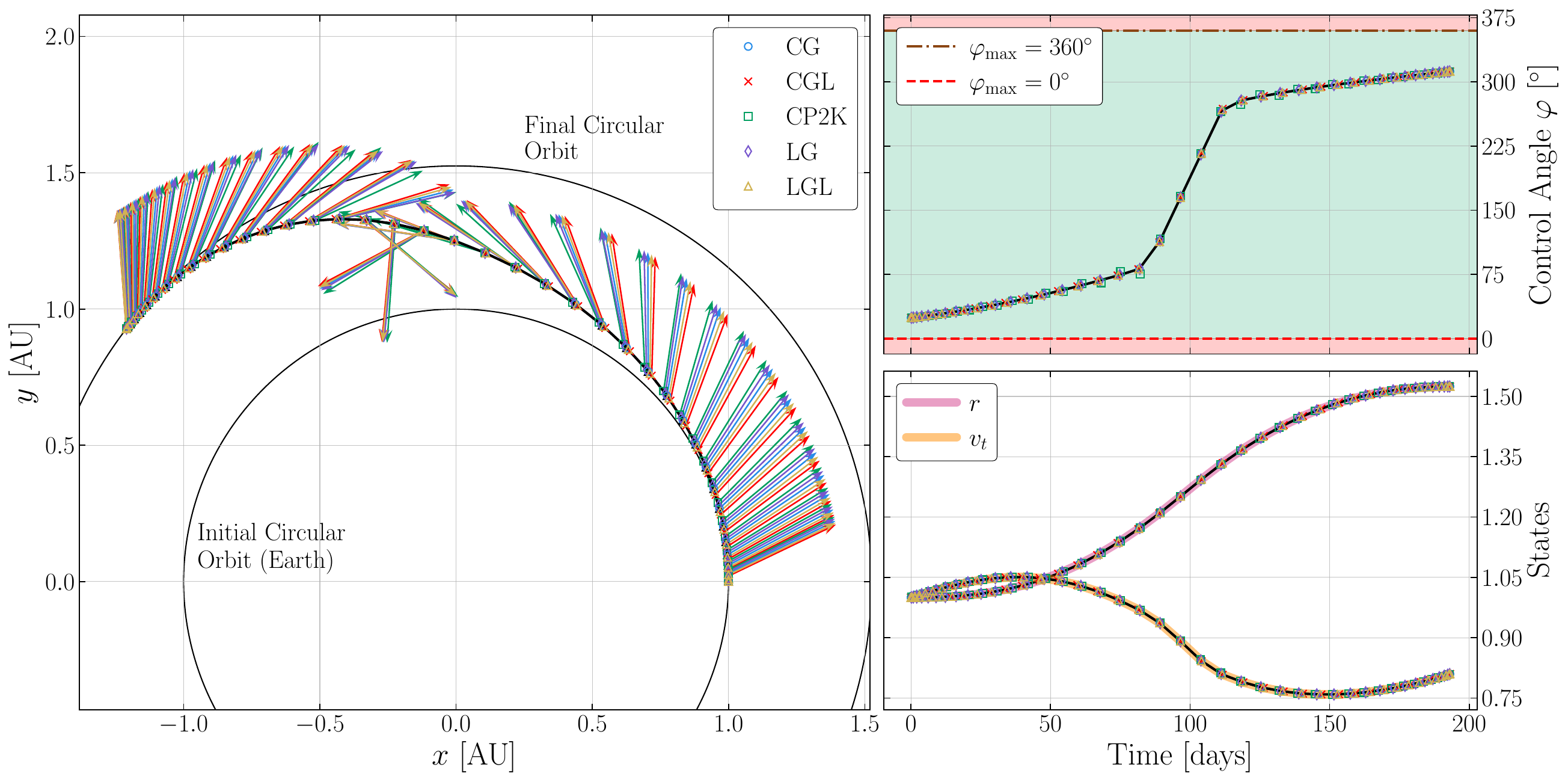}
    \caption{Maximum radius planar circular orbit raising under constant thrust solved with $n=40$. Constrained regions are shown in red, while unconstrained region is shown in green.}
    \label{fig:OrbitRaisingMaxRad}
\end{figure}

As with the minimum time problem, the maximum radius problem is broken up into two sub-problems, and all parameters were initialized in the same manner. Again, 100 trials were run for $n=40$. The number of successful runs, average numbers of iterations, final orbit radius, and computational runtimes are given in Table~\ref{tab:MaxRadius}. This problem appears to be more sensitive to the polynomial and node choice, as exhibited by the difference in the number of successful runs between the CG/CGL/CP2K nodes and LG/LGL nodes. The phase space and control and state histories for a converged solution are shown in Fig.~\ref{fig:OrbitRaisingMaxRad}, and as expected, the results agree with those of the minimum time problem. The worse convergence of this problem compared to the minimum time problem is believed to be due to the different cost function and constraints. Without the ability to alter the final time, the coefficients have less freedom to satisfy the terminal constraints, which, instead of being known values, are now functions of the final state of the trajectory.

\begin{table}[htb!]
    \caption{\label{tab:MaxRadius} DOPIC Maximum Radius Orbit Transfer Performance Results Averaged Over 100 Runs for $n=40$}
    \centering
    \begin{tabular}{lcccc}
        \hline \hline
        \mc{\textbf{Nodes}} & \mc{\textbf{Successful Runs}} & \mc{\textbf{Iterations}} & \mc{\textbf{Final Radius $r_f$ [ND]}} & \mc{\textbf{Runtime [s]}} \\ \hline
        CG   & 68  & 246 & 1.5246545 & 1.667 \\
        CGL  & 69  & 337 & 1.5246583 & 1.974 \\
        CP2K & 63  & 248 & 1.5248996 & 1.614 \\
        LG   & 94  & 296 & 1.5246556 & 1.815 \\
        LGL  & 97  & 378 & 1.5246579 & 2.137 \\
        \hline \hline
    \end{tabular}
\end{table}

\subsection{Rocket Landing Flip Maneuver}
In the last decade, there have been significant technological advances which have allowed for landing, recovery, and reuse of orbital class rockets. SpaceX has successfully launched and landed over 500 Falcon 9 vehicles, both on land at Cape Canaveral and on autonomous floating drone ships in the Atlantic and Pacific Oceans \cite{SXF9}. In order to enable precision landing, the vehicle must continuously recompute an optimal trajectory to ensure that it can reach the landing pad while accounting for unmodeled dynamics and potential system anomalies. To accomplish this, convex optimization techniques are used to guarantee the optimal solution can be found quickly \cite{BlackmoreAut, CVXGEN}.

SpaceX is currently developing the Starship vehicle to be the first fully and rapidly reusable rocket which will return humans back to the Moon and eventually to Mars \cite{SXSS, ShotwellBlackmore}. Following a successful launch, the upper stage (also referred to as Starship) will perform a controlled atmospheric reentry using the forward and aft flaps, and initiate a landing burn to reorient the vehicle and land vertically at the desired location. This landing flip maneuver requires relighting the sea-level Raptor 3 engine(s), and has been successfully demonstrated \cite{SN15}. Several of the same convexification techniques used on Falcon 9 have been applied to the Starship landing problem \cite{DualQuat, RealTime}. Here, a simplified model with a single Raptor 3 engine is considered without convexification. Extensive runtime metrics are not reported on this problem, as its purpose is to demonstrate OPIC's ability to solve very complex, constrained OCPs.

The dynamics for this problem are
\begin{subequations} \label{eq:ddots}
\begin{align}
    \ddot{h}(t) &= \frac{\delta \Tmax \cos(\theta + \varphi) - D(t) \dot{h}}{m} - g_0 \\
    \ddot{x}(t) &= \frac{-\delta \Tmax \sin(\theta + \varphi) - D(t) \dot{x}}{m} \\
    \ddot{\theta}(t) &= \frac{-\lcm \delta \Tmax \sin\varphi + M_A(t) }{I(t)} \\
    \dot{m}(t) &= - \frac{\delta \Tmax}{\Isp g_0}
\end{align}
\end{subequations}
The vehicle states are the altitude $h$, downrange position $x$, attitude $\theta$ (measured from the vertical inertial frame), and mass $m$. The engine parameters are the throttle $\delta$, maximum thrust $T_{\max}$, nozzle gimbal angle $\varphi$, and specific impulse $\Isp$. The drag function $D(t)$, aerodynamic moment $M_A(t)$ and rotational inertia $I(t)$ defined in Eqs.~\eqref{eq:Drag}-\eqref{eq:Moment} in terms of the vehicle length $\lh$, radius $\lr$, center of mass $\lcm$, center of pressure $\lcp$, and a single aerodynamic coefficient $C_{L/D}$. The atmospheric density and gravitational acceleration at the Earth's surface are given by $\rho_0$ and $g_0$, respectively. The control variables are the engine throttle $\delta$ and gimbal angle $\varphi$.
\begin{align}
    D(t) &= \frac{1}{2} \rho_0 \aref C_{L/D} \sqrt{\dot{h}^2 + \dot{x}^2} \label{eq:Drag} \\
    I(t) &= m \left( \frac{\lr^2}{4} + \frac{\lh^2}{12} \right) \label{eq:MOI} \\
    M_A(t) &= (\lcp - \lcm) D(t) (\dot{x} \cos\theta + \dot{h}\sin\theta) \label{eq:Moment}
\end{align}
The optimal control problem is then defined as the free final time, minimum fuel problem,
\begin{align}
    \min \quad J &= -m(t_f) \\
    \begin{split}
        \text{s.t.} \quad &\text{Eqs.~\eqref{eq:ddots}} \\
        h(0) &= h_0, \quad \dot{h}(0) = \dot{h}_0, \quad x(0) = x_0, \quad \dot{x}(0) = \dot{x}_0, \quad \theta(0) = \theta_0, \quad \dot{\theta}(0) = \dot{\theta}_0, \quad m(0) = m_0 \\
        h(t_f) &= h_f, \quad \dot{h}(t_f) = \dot{h}_f, \quad x(t_f) = x_f, \quad \dot{x}(t_f) = \dot{x}_f, \quad \theta(t_f) = \theta_f, \quad \dot{\theta}(t_f) = \dot{\theta}_f \\
        \delta_{\min} &\leq \delta \leq \delta_{\max} \\
        |\varphi| &\leq \varphi_{\max} \\
        |\dot{\varphi}| &\leq \dot{\varphi}_{\max} \\
        m(t_f) &\geq m_\text{dry}
    \end{split}
\end{align}
where $\delta_{\min}$ and $\delta_{\max}$ are the minimum and maximum engine throttles, $\varphi_{\max}$ is the maximum nozzle gimbal angle, $\dot{\varphi}_{\max}$ is the maximum nozzle gimbal rate, and $m_\text{dry}$ is the vehicle's dry mass.

Given the large difference in scales between various states and parameters, the problem is non-dimensionalized with the length-scale $L := h_0 \: [\unit[]{\meter}]$, mass-scale $M := m_0 \: [\unit[]{\kilogram}]$, and time-scale $\beta := \sqrt{h_0 / g_0} \: [\unit[]{\second}]$. Thus, denoting the scaled states with an overbar $\bar{(\cdot)}$, the dimensionless ODEs are
\begin{subequations} \label{eq:ddots_bar}
\begin{align}
    \ddot{\bar{h}} &= c_1 \left( \frac{\delta \cos(\bar{\theta}+\bar{\varphi})}{\bar{m}} \right) - c_2 \left( \frac{\sqrt{\dot{\bar{h}}^2 + \dot{\bar{x}}^2} \dot{\bar{h}}}{\bar{m}} \right) - c_3 \label{eq:hddotbarfinal} \\
    \ddot{x} &= -c_1 \left( \frac{\delta \sin(\bar{\theta}+\bar{\varphi})}{\bar{m}} \right) - c_2 \left( \frac{\sqrt{\dot{\bar{h}}^2 + \dot{\bar{x}}^2} \dot{\bar{x}}}{\bar{m}} \right) \label{eq:xddotbarfinal} \\
    \ddot{\bar{\theta}} &= -c_4 \left( \frac{\delta \sin\bar{\varphi}}{\bar{m}} \right) + c_5 \left( \frac{\sqrt{\dot{\bar{h}}^2 + \dot{\bar{x}}^2} (\dot{\bar{x}} \cos\bar{\theta} + \dot{\bar{h}} \sin\bar{\theta})}{\bar{m}} \right) \label{eq:thetaddotbarfinal} \\
    \dot{\bar{m}} &= -c_6 \delta \label{eq:mdotbarfinal}
\end{align}
\end{subequations}
where the constants are defined as
\begin{equation}
\begin{aligned}
    c_1 &:= \frac{\beta^2 \Tmax}{L M}, \qquad &&c_2 := \frac{\rho_0 \aref C_{L/D} L}{2 M}, \qquad &&c_3 := \frac{\beta^2 g_0}{L} \\
    c_4 &:= \frac{\beta^2 \lcm \Tmax}{M \left( \frac{\lr^2}{4} + \frac{\lh^2}{12} \right)}, \qquad &&c_5 := \frac{(\lcp - \lcm) \rho_0 \aref C_{L/D} L^2}{2 M \left( \frac{\lr^2}{4} + \frac{\lh^2}{12} \right)}, \qquad &&c_6 := \frac{\beta \Tmax}{M \Isp g_0}
\end{aligned}
\end{equation}
The non-dimensionalized dynamics are then approximated at the node set $\bm{\tau}$ in the computational domain using OPIC. The final OCP is then written as
\begin{align}
    \min \quad J &= -\bar{m}(t_f) \\
    \text{s.t.} \quad & \text{Eqs.~\eqref{eq:ddots_bar}}
\end{align}
The boundary conditions and inequality constraints, not shown here for brevity, are simply transformed in the same manner as the states and controls. All of the parameters used for this problem, gathered or estimated from various public sources, are given in Table~\ref{tab:RocketEnvironParams}. 

\begin{table}[htb!]
    \caption{\label{tab:RocketEnvironParams} Atmospheric Constants, Starship Vehicle Parameters, and Boundary Condition Values \cite{SXSS, RealTime}}
    \centering
    \begin{tabular}{llr}
        \hline \hline
        \ml{\textbf{Variable Name}} & \ml{\textbf{Symbol}} & \mr{\textbf{Value}} \\ \hline
        Atmospheric density & $\rho_0$ & $1.225~\unit[per-mode=symbol]{\kilogram\per\meter\cubed}$ \\
        Gravitational acceleration & $g_0$ & $9.81~\unit[per-mode=symbol]{\meter\per\second\squared}$ \\
        Dry mass & $m_\text{dry}$ & $85,000~\unit{\kilogram}$ \\
        Vehicle length & $\lh$ & $50~\unit{\meter}$ \\
        Vehicle radius & $\lr$ & $4.5~\unit{\meter}$ \\
        Center of mass & $\lcm$ & $20~\unit{\meter}$ \\
        Center of pressure & $\lcp$ & $22.5~\unit{\meter}$ \\
        Aerodynamic coefficient & $C_{L/D}$ & $0.4$ \\
        Specific impulse & $I_\text{sp}$ & $350~\unit{\second}$ \\
        Max thrust & $T_{\max}$ & $280~\text{tf}$ \\
        Max throttle & $\delta_{\max}$ & $1$ \\
        Min throttle & $\delta_{\min}$ & $0.4$ \\
        Max nozzle gimbal angle & $\varphi_{\max}$ & $15~\unit{\deg}$ \\
        Max nozzle gimbal rate & $\dot{\varphi}_{\max}$ & $15~\unit[per-mode=symbol]{\deg\per\second}$ \\
        Initial altitude & $h_0$ & $1000~\unit{\meter}$ \\
        Initial vertical velocity & $\dot{h}_0$ & $-90~\unit[per-mode=symbol]{\meter\per\second}$ \\
        Initial downrange & $x_0$ & $100~\unit{\meter}$ \\
        Initial downrange velocity & $\dot{x}_0$ & $0~\unit[per-mode=symbol]{\meter\per\second}$ \\
        Initial attitude & $\theta_0$ & $\sfrac{\pi}{2}~\unit{\rad}$ \\
        Initial attitude rate & $\dot{\theta}_0$ & $0~\unit[per-mode=symbol]{\rad\per\second}$ \\
        Initial mass & $m_0$ & $100,000~\unit{\kilogram}$ \\
        Final altitude & $h_f$ & $0~\unit{\meter}$ \\
        Final vertical velocity & $\dot{h}_f$ & $0~\unit[per-mode=symbol]{\meter\per\second}$ \\
        Final downrange & $x_f$ & $0~\unit{\meter}$ \\
        Final downrange velocity & $\dot{x}_f$ & $0~\unit[per-mode=symbol]{\meter\per\second}$ \\
        Final attitude & $\theta_f$ & $0~\unit{\rad}$ \\
        Final attitude rate & $\dot{\theta}_f$ & $0~\unit[per-mode=symbol]{\rad\per\second}$ \\
        \hline \hline
    \end{tabular}
\end{table}

To solve this problem using DOPIC, fourth $n\tth$ order polynomials series are used to approximate the non-dimensionalized dynamics, with engine throttle and nozzle gimbal angle controls at the interpolation nodes. The first three series represent the second order dynamics $\ddot{\bar{h}}$, $\ddot{\bar{x}}$, and $\ddot{\bar{\theta}}$, and the last series represents the first order mass dynamics $\dot{\bar{m}}$. The constraints on the nozzle gimbal rate $\bm{c}_{\dot{\bar{\varphi}}_{\max}}$ are defined in the same way as the previous problem with Eq.~\eqref{eq:att_rate_con}. Thus, the free variable vector $\chi$, equality constraint vector $\Ceq$, and inequality constraint vector $\Cin$ are given as
\begin{align}
    \chi &= \begin{bmatrix} \bm{\alpha}_{\bar{h}} \\ \bm{\alpha}_{\bar{x}} \\ \bm{\alpha}_{\bar{\theta}} \\ \bm{\alpha}_{\bar{m}} \\ \bm{\bar{\varphi}} \\ \bm{\delta} \\ t_f \end{bmatrix} \in \mathbb{R}^{(6(n+1)+1) \times 1}, \quad \Ceq = \begin{bmatrix} \bar{h}(-1) - \bar{h}_0 \\ \bar{h}^{\prime}(-1) - \left( \frac{t_f}{2} \right) \dot{\bar{h}}_0 \\ \bar{x}(-1) - \bar{x}_0 \\ \bar{x}^{\prime}(-1) - \left( \frac{t_f}{2} \right) \dot{\bar{x}}_0 \\ \bar{\theta}(-1) - \bar{\theta}_0 \\ \bar{\theta}^{\prime}(-1) - \left( \frac{t_f}{2} \right) \dot{\bar{\theta}}_0 \\ \bar{m}(-1) - \bar{m}_0 \\ \Xi_{\bar{h}} \\\Xi_{\bar{x}} \\ \Xi_{\bar{\theta}} \\ \Xi_{\bar{m}} \\ \bar{h}(1) - \bar{h}_f \\ \bar{h}^{\prime}(1) - \left( \frac{t_f}{2} \right) \dot{\bar{h}}_f \\ \bar{x}(1) - \bar{x}_f \\ \bar{x}^{\prime}(1) - \left( \frac{t_f}{2} \right) \dot{\bar{x}}_f \\ \bar{\theta}(1) - \bar{\theta}_f \\ \bar{\theta}^{\prime}(1) - \left( \frac{t_f}{2} \right) \dot{\bar{\theta}}_f \end{bmatrix} = \bm{0} \in \mathbb{R}^{(4(n+1)+13) \times 1}, \quad \Cin = \begin{bmatrix} \bm{\delta} - \delta_{\max} \\ \delta_{\min} - \bm{\delta} \\ \bm{\bar{\varphi}} - \bar{\varphi}_{\max} \\ -\bar{\varphi}_{\max} - \bm{\bar{\varphi}} \\ \bm{\bar{m}} - \bar{m}_0 \\ \bar{m}_\text{dry} - \bm{\bar{m}} \\ \bm{c}_{\dot{\bar{\varphi}}_{\max}} \\ -t_f \end{bmatrix} \leq \bm{0} \in \mathbb{R}^{7(n+1) \times 1}
\end{align}

\begin{figure}[H]
    \centering
    \includegraphics[width=.95\textwidth]{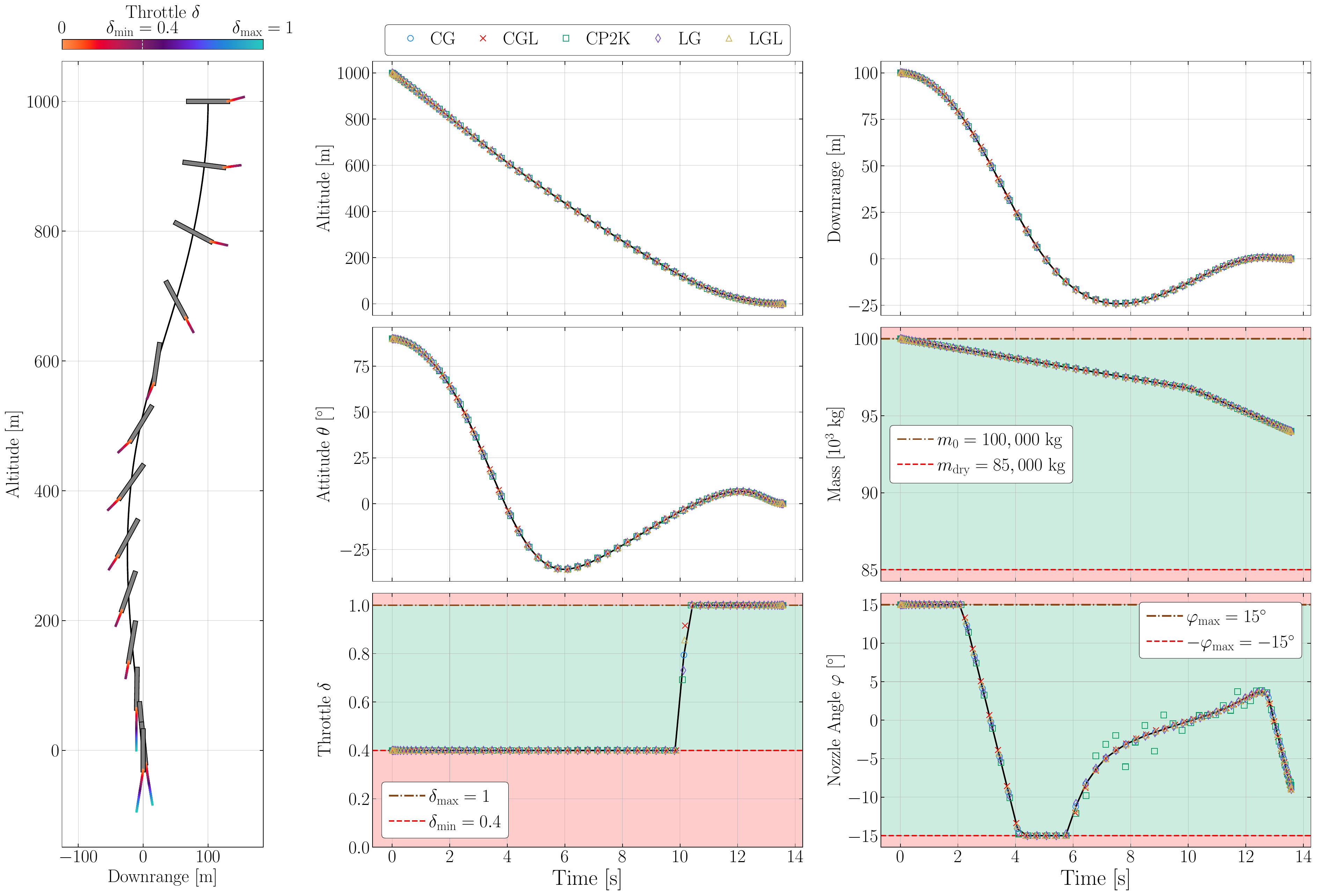}
    \caption{Phase space and state and control histories for the rocket landing flip maneuver solved with $n=60$. Constrained regions are shown in red, while unconstrained regions are shown in green.}
    \label{fig:RocketLanding}
\end{figure}

The non-dimensionalized problem is solved using an order $n = 60$ approximation, and is then transformed back to the dimensional states. The initial guess for the free variable vector is as follows: all coefficients are randomly distributed between zero and one, the gimbal angles are zero, the throttles are $0.8$, and the final time is 10 seconds. The Starship vehicle trajectory are shown alongside the state and control histories and constraints in Fig.~\ref{fig:RocketLanding}. The performance index (final mass $m_f$) and final time $t_f$ for each node set are given in Table~\ref{tab:StarshipPerf}. All solutions are similar, although some oscillations are seen in the nozzle gimbal angle for CP2K. This is a result of the inexact gimbal angle rate approximation, and can be alleviated with a higher order approximation or reparameterization of the control angle. The final mass and time for all nodes is seen to be $m_f \simeq 93979$~kg and $t_f \simeq 13.57$~s, respectively.

\begin{table}[H]
    \caption{\label{tab:StarshipPerf} DOPIC Starship Landing Flip Maneuver Performance Results $n=60$}
    \centering
    \begin{tabular}{lcc}
        \hline \hline
        \mc{\textbf{Nodes}} & \mc{\textbf{Final Mass} $m_f~[\unit{\kilogram}]$} & \mc{\textbf{Final Time} $t_f~[\unit{\second}]$} \\ \hline
        CG   & 93979.49 & 13.5794 \\
        CGL  & 93979.20 & 13.5668 \\
        CP2K & 93979.86 & 13.5804 \\
        LG   & 93979.98 & 13.5861 \\
        LGL  & 93979.22 & 13.5730 \\
        \hline \hline
    \end{tabular}
\end{table}


%% file: 05_conclusion.tex
\section{Conclusion} \label{sec:conclusion}
This paper presents a generalized direct pseudospectral method to transcribe the continuous time optimal control problem based on integration rather than differentiation. Orthogonal polynomial integral collocation allows for the approximation of the highest order derivative so that only a single set of polynomial coefficients is necessary to represent any state derivative level. This greatly reduces the parameter space when the problem is solved using a nonlinear programming algorithm. The proposed method can be used with any orthogonal polynomial set and associated grid. Example problems are solved using Chebyshev polynomials of the first kind, Chebyshev polynomials of the second kind, and Legendre polynomials with common grids associated with their zeros and/or extrema. Gaussian quadrature is employed to approximate integral valued cost functionals.

Example problems with state and control inequality constraints that exhibit analytical solutions are solved to verify the proposed solution method. Chebyshev polynomials of the second kind are shown to outperform the other polynomials tested in average runtime. Next, the planar circular orbit raising problem is solved for both the minimum time and maximum radius cases. Solutions using all polynomials are shown to agree with other pseudospectral methods, though Legendre polynomials and grids displayed significantly higher convergence rates when compared to Chebyshev polynomials. Lastly, a larger dimensional, highly constrained problem is solved in the form of a minimum fuel mass, atmospheric rocket flip landing maneuver to demonstrate the proposed method's capability. Results demonstrate the efficacy and efficiency with which the proposed method can solve problems that contain multiple high-order derivatives because of the integral approach's reduction of parameter space.